\documentclass[11pt]{article}
\usepackage[T1]{fontenc}

\usepackage[utf8]{inputenc}

\usepackage{mathrsfs}
\usepackage{amsmath}
\usepackage{amssymb}
\usepackage{graphicx}
\usepackage{color}
\usepackage{caption,subcaption}
\usepackage{babel}
\usepackage{comment}
\usepackage{algorithm}
\usepackage{algpseudocode}
\usepackage{setspace}
\usepackage{grffile}
\usepackage{enumerate}

\usepackage{hyperref}
\hypersetup{colorlinks=true,
            linkcolor=black,
            anchorcolor=black,
            citecolor=black}

\numberwithin{equation}{section}
\numberwithin{figure}{section}
\numberwithin{table}{section}

\newcommand{\R}{\mathbb{R}}

\newcommand{\lam}{\lambda}

\newcommand{\bx}{\mathbf{x}}

\newcommand{\gam}{\gamma}
\newcommand{\bv}{\mathbf{v}}
\renewcommand{\th}{\theta}

\newcommand{\eps}{\varepsilon}

\newcommand{\network}{\mathcal{N}}

\newcommand{\bE}{\mathbb{E}}
\newcommand{\bB}{\mathbb{B}}
\newcommand{\bH}{\mathbb{H}}

\newcommand{\Dx}{\Delta x}
\newcommand{\Dt}{\Delta t}

\newcommand{\iph}{i+\frac{1}{2}}
\newcommand{\imh}{i-\frac{1}{2}}

\newcommand{\nph}{n+\frac{1}{2}}
\newcommand{\nmh}{n-\frac{1}{2}}
\newcommand{\rd}{\mathrm{d}}

\newcommand{\hE}{\widehat{E}}
\newcommand{\hrho}{\hat{\rho}}
\newcommand{\ion}{\rho_{ion}}
\newcommand{\hphi}{\widehat{\phi}}

\newcommand{\hH}{\widehat{H}}
\newcommand{\hB}{\widehat{B}}
\newcommand{\hh}{\widehat{h}}

\newcommand{\cJ}{\mathcal{J}}
\newcommand{\cL}{\mathcal{L}}

\newcommand{\feq}{f_{eq}}

\newcommand{\Beq}{B_{eq}}
\newcommand{\bEeq}{\mathbb{E}_{eq}}
\newcommand{\bBeq}{\mathbb{B}_{eq}}

\newcommand{\energy}{\mathcal{E}}

\newcommand{\dE}{\delta E}
\newcommand{\dB}{\delta B}

\newcommand{\df}{\delta f}

\newcommand{\drho}{\delta \rho}
\newcommand{\dH}{\delta H}
\newcommand{\laplace}{\mathscr{L}}

\newcommand{\baB}{\overline{B}}
\newcommand{\Bext}{B_{\mathrm{ext}}}
\newcommand{\bBext}{\bB_{\mathrm{ext}}}
\newcommand{\dJ}{\delta J}
\newcommand{\bav}{\overline{v}}

\newcommand{\ds}{\displaystyle}

\newcommand{\email}[1]{\href{mailto:#1}{\nolinkurl{#1}}}

\title{Structure-Informed Electromagnetic Feedback Control of the Vlasov–Maxwell System with Delayed Observations}

\author{Jingcheng Lu \thanks{Institute of Mathematical Sciences, ShanghaiTech University, China; School of Mathematics, University of Minnesota–Twin Cities, MN (\email{lujch@shanghaitech.edu.cn}).}
\and Li Wang \thanks{School of Mathematics, University of Minnesota–Twin Cities, MN (\email{liwang@umn.edu}).}
}
\date{}

\begin{document}
\maketitle

\begin{abstract}
   Feedback control of kinetic plasma instabilities is challenging because observing the full phase-space distribution is costly, while sensing and actuation inevitably introduce delay. We develop a structure-informed electromagnetic feedback framework for the Vlasov--Maxwell system that uses only macroscopic observations to control the full kinetic distribution. We first show that magnetic actuation alone does not directly regulate radial velocity moments, motivating the use of combined electric and magnetic controls. A nonlinear perturbation estimate yields an instantaneous cancellation--dissipation feedback law. To account for delayed observations, we then use Maxwell's equations to identify the macroscopic quantities governing the resulting cancellation errors, namely currents and spatial field gradients. These quantities serve as inputs to a low-rank neural operator that learns a residual correction to the delayed cancellation baseline. Thus, while the feedback depends only on low-dimensional macroscopic observables, the control objective remains the stabilization of the full distribution toward a prescribed kinetic equilibrium. Numerical experiments for two-stream and Weibel-type instabilities demonstrate effective long-time stabilization under substantial feedback delays.
\end{abstract}

\section{Introduction}
Controlling plasma instabilities is an important challenge in applications
ranging from fusion energy to charged-particle beam manipulation. External
electromagnetic fields provide a means of influencing the plasma, but their
effectiveness depends on the dynamics and the control objective.
Magnetic confinement and macroscopic stability have long been central
topics in fusion research \cite{zohm2014magnetohydrodynamic}. At the kinetic
level, the Vlasov--Poisson (VP) and Vlasov--Maxwell (VM) equations describe
the interaction between the particle distribution and self-consistent
fields. Instabilities such as the two-stream, bump-on-tail, and Weibel
instabilities illustrate how small perturbations of an unstable
distribution can generate growing fields and substantial changes in the
plasma state. Their suppression calls for control strategies that account
for this coupled evolution.

The mathematical literature on kinetic plasma control includes both
controllability and optimal-control formulations. Controllability results
for the VP system using localized interior source controls were
established in \cite{glass2003controllability} and extended to systems
with prescribed external force fields in \cite{glass2012external}.
For electromagnetic actuation, a variational framework for optimal
magnetic control of a VP plasma was developed in
\cite{knopf2018optimal}, with an extension to magnetic fields generated
by a finite collection of fixed coils in \cite{knopf2020coils}.
Adjoint-based optimization has also been combined with Monte Carlo
particle-in-cell discretizations for magnetic confinement
\cite{bartsch2024controlling}. Optimal control has been studied for
the two-dimensional relativistic VM system, with external currents acting
as controls in Maxwell's equations \cite{weber2021optimal}.

For kinetic instability suppression, a PDE-constrained optimization
framework was proposed in \cite{einkemmer2024suppressing}, in which a
parameterized time-independent external electric field is chosen to
maintain a desired distribution and delay the onset of VP instabilities.
To reduce the computational cost of repeated kinetic solves, a
moment-based formulation uses a reduced moment system during
optimization, with the resulting controls assessed against the full
kinetic dynamics \cite{lu2025controlling}. The choice of objective is
also consequential: numerical investigations of the associated
optimization landscapes indicate that objectives incorporating
time-integrated information can be more favorable for gradient-based
optimization \cite{guerra2026metric}. These studies highlight the importance
of both the forward model and the objective in constructing effective
controls.

A complementary approach uses linear stability analysis to identify
and suppress unstable modes. Fourier--Laplace analysis of the VP system
leads to electric-field controls designed to cancel unstable contributions
in the linear response \cite{einkemmer2024control}. This perspective
has been extended to a plasma subject to a uniform background magnetic
field, with electric controls designed through the dispersion relation
\cite{chen2026control}. For electromagnetic kinetics, control by
prescribed transverse laser fields has been studied in a reduced VM
model \cite{crouseilles2025control}. This approach combines
linear analysis for transverse instabilities with nonlinear
PDE-constrained optimization for longitudinal instability suppression.
In particular, the nonlinear coupling can be exploited when the
corresponding longitudinal mode is unaffected by the laser control
in the linearized system. The control fields in these constructions
are determined from the specified initial data and subsequently applied
without updating them from observations of the evolving plasma.

Feedback control provides a way to incorporate such observations into
the control action. At the kinetic level, instantaneous magnetic feedback
laws have been derived from short-time optimization of particle dynamics
to maintain spatial confinement \cite{albi2025instantaneous}.
This approach was extended to
collisional VP--BGK dynamics with uncertain parameters in
\cite{albi2025robust}. Learning-based feedback has also been demonstrated
in tokamak experiments: deep reinforcement learning has been used for
magnetic control of plasma shape and position
\cite{degrave2022magnetic} and for avoidance of tearing instabilities
\cite{seo2024avoiding}. These developments motivate the construction
of feedback policies that can be evaluated rapidly while adapting to
the observed state.

More directly related to the present work, an operator-based feedback
framework for the VP system was developed in \cite{lu2025dynamical}.
The framework uses a low-rank neural operator that maps the distribution
perturbation to an external electric field. The operator is trained through an adjoint
method and demonstrates effective stabilization numerically beyond the
optimization horizon. An analytical cancellation-based control law
guided by a perturbation energy estimate was also introduced in that
work. This operator viewpoint
is consistent with the broader development of neural operators for
maps between function spaces, including Fourier neural operators,
DeepONet, and integral-kernel architectures
\cite{li2020fourier,lu2021learning,kovachki2023neural}. Extending the VP
feedback framework to electromagnetic kinetics raises additional
questions about the respective roles of electric and magnetic actuation,
the information required by the controller, and the effect of feedback
latency.


To be more specific, we consider the nonrelativistic VM system. In normalized variables,
the uncontrolled system reads
\begin{equation}\label{eq:VM-3D3V}
\left\{
\begin{aligned}
&\partial_t f+\bv\cdot\nabla_{\bx}f
 +(\bE+\bv\times\bB)\cdot\nabla_{\bv}f=0,\\
&\partial_t\bE=\nabla_{\bx}\times\bB-\mathbf J,\\
&\partial_t\bB=-\nabla_{\bx}\times\bE,\\
&\nabla_{\bx}\cdot\bE=\rho-\ion,\qquad
 \nabla_{\bx}\cdot\bB=0.
\end{aligned}
\right.
\end{equation}
Here $f=f(\bx,\bv,t)$ is the particle distribution,
$\bx\in\Omega=\mathbb T^3$ and $\bv\in\R^3$ denote position and
velocity, and $\bE(\bx,t)$ and $\bB(\bx,t)$ are the self-consistent
electric and magnetic fields. The density and current are
$\rho(\bx,t)=\int_{\R^3}f(\bx,\bv,t)\,d\bv$ and
$\mathbf J(\bx,t)=\int_{\R^3}\bv f(\bx,\bv,t)\,d\bv$,
respectively, and $\ion$ is the fixed neutralizing background.
We introduce external electric and magnetic fields in the Lorentz
force and seek feedback laws that suppress the growth of perturbations
around a prescribed spatially homogeneous equilibrium.

Two issues guide the controller design. First, electric and magnetic actuation have different effects on the kinetic state. The magnetic Lorentz force performs no work and makes no direct contribution to the evolution of radial velocity moments. Although magnetic control can influence these moments indirectly through the self-consistent fields, it does not provide a general mechanism for regulating an arbitrary kinetic-energy profile. This motivates combining electric and magnetic control channels. Second, the controller must operate with {\it limited} and {\it delayed} observations. The learned feedback controller in \cite{lu2025dynamical} takes the full phase-space perturbation $f-f_{\mathrm{eq}}$ as its input. In contrast, the feedback law developed here depends only on macroscopic observables, including current densities, electromagnetic fields, and spatial gradients of the transverse fields. Importantly, these macroscopic quantities are used only as feedback information: the control objective remains kinetic, namely, to suppress the perturbation of the full distribution $f$ and drive it toward the prescribed equilibrium $f_{\mathrm{eq}}$. Thus, the proposed framework seeks to {\it control a high-dimensional kinetic state using only low-dimensional macroscopic observations, without requiring direct observation or reconstruction of the full phase-space distribution during deployment}. Moreover, sensing, computation, and actuation introduce latency, which can undermine the effectiveness of instantaneous field-cancellation strategies when they are implemented using delayed observations.

We address these issues through a structure-informed feedback design.
The nonlinear $L^2$ perturbation balance provides an instantaneous
cancellation--dissipation prototype. Maxwell's equations then identify
the currents and transverse field gradients that generate the
cancellation errors under delayed feedback. We use these quantities
as inputs to a low-rank linear neural operator that learns a residual
correction to the delayed field-cancellation baseline. The deployed
controller uses current and field observations, while its offline
training uses kinetic simulations and an adjoint-based optimization
method. This construction connects the choice of feedback inputs
and the learned correction to the governing dynamics.

The main contributions are as follows:
\begin{itemize}

\item We identify a structural limitation of purely magnetic actuation for kinetic-state stabilization. In particular, magnetic fields act through velocity-space rotations and do not directly modify radial velocity moments of order $k\geq 2$, including the kinetic energy. This shows that magnetic fields' ability to alter particle trajectories does not, in general, provide sufficient control authority to drive the distribution toward a prescribed equilibrium, motivating the use of hybrid electromagnetic feedback.

\item Based on the Vlasov--Maxwell energy structure, we construct a cancellation--dissipation feedback framework that combines electric and magnetic control channels. To accommodate delayed observations, we further derive a structure-informed compensation strategy from Maxwell's equations, using current and field-gradient observables as feedback inputs and a low-rank neural operator to learn the corresponding residual correction.

\item Numerical experiments demonstrate that the proposed feedback remains effective well beyond the optimization horizon and under substantial feedback delays. Comparisons with alternative observable choices further show that incorporating the governing Vlasov--Maxwell structure into the controller design is essential for long-time stabilization.

\end{itemize}

\section{Motivation for electromagnetic feedback control}

A fundamental question in plasma feedback control is the choice of
actuation. Magnetic fields play a central role in plasma confinement
and equilibrium shaping \cite{zohm2014magnetohydrodynamic,albi2025instantaneous}:
they bend particle trajectories without doing work on individual
particles. Here, our objective is to suppress the growth of kinetic
perturbations around a prescribed equilibrium. This raises the question
of whether purely magnetic actuation provides sufficient control
authority for general kinetic-state regulation.

To answer this question, we denote the external electric
and magnetic controls by $\bH(\bx,t)$ and $\bBext(\bx,t)$, then the controlled Vlasov equation is
\begin{equation}\label{eq:VM-control-3D3V}
\partial_t f+\bv\cdot\nabla_{\bx}f
 +\big[\bE+\bH+\bv\times(\bB+\bBext)\big]
   \cdot\nabla_{\bv}f=0.
\end{equation}
We assume smooth solutions, periodic spatial boundary conditions,
and sufficient decay in velocity to justify the following
integrations by parts and moment calculations.

Let $(\feq,\bEeq,\bBeq)$ be a stationary equilibrium of the
uncontrolled system, and define
\begin{equation*}
\df=f-\feq,\qquad
\delta\bE=\bE-\bEeq,\qquad
\delta\bB=\bB-\bBeq.
\end{equation*}
Subtracting the equilibrium equation from
\eqref{eq:VM-control-3D3V} gives
\begin{equation}\label{eq: VM perturbation}
\begin{split}
&\partial_t\df+\bv\cdot\nabla_{\bx}\df
 +\big[\bE+\bH+\bv\times(\bB+\bBext)\big]
   \cdot\nabla_{\bv}\df\\
&\qquad
 +\big[\delta\bE+\bH+\bv\times(\delta\bB+\bBext)\big]
   \cdot\nabla_{\bv}\feq=0.
\end{split}
\end{equation}
For every velocity-independent vector field $\mathbf C(\bx,t)$,
\begin{equation*}
\bv\cdot(\bv\times\mathbf C)=0,\qquad
\nabla_{\bv}\cdot(\bv\times\mathbf C)=0.
\end{equation*}
In particular, the phase-space transport is divergence free.
Multiplying \eqref{eq: VM perturbation} by $\df$ and integrating
over $\Omega\times\R^3$ yields the exact nonlinear perturbation balance
\begin{equation}\label{eq:energy estimate}
\begin{split}
\frac{1}{2}\frac{d}{dt}\|\delta f(t)\|_{L^2}^2
=& -\int (\delta\bE+\bH)\cdot
      \left[\int \delta f\,\nabla_{\bv}\feq\, d\bv\right]d\bx
  +\int (\delta\bB+\bB_{\mathrm{ext}})\cdot
      \left[\int (\bv\times\nabla_{\bv}\feq)\,\delta f\, d\bv\right]d\bx\\
\end{split}
\end{equation}
This identity measures the evolution of the $L^2$ perturbation.
Electric actuation couples to the velocity gradient of the equilibrium,
while magnetic actuation couples to its angular derivatives through
$\bv\times\nabla_{\bv}\feq$. The latter vanishes when the equilibrium
is isotropic in velocity.

The distinction between electric work and magnetic rotation is also
visible in the radial velocity moments
\begin{equation*}
M_k(t)=\int_\Omega\int_{\R^3}|\bv|^k f(\bx,\bv,t)\,d\bv\,d\bx,
\qquad k\geq2.
\end{equation*}
Using the orthogonality and divergence identities above, integration
by parts gives
\begin{equation}\label{eq:radial-moment-3D3V}
\frac{\rd}{\rd t }M_k(t)
=k\int_\Omega\int_{\R^3}|\bv|^{k-2}
  \bv\cdot(\bE+\bH)f\,d\bv\,d\bx.
\end{equation}
Thus, the magnetic force makes no direct contribution to these
moments, including the kinetic energy $M_2/2$. Magnetic actuation
can nevertheless affect them indirectly by changing the distribution,
current, and self-consistent electric field.

A simple example illustrates the resulting limitation on general
kinetic-state regulation. Consider a neutral, spatially homogeneous,
isotropic initial distribution with zero initial self-consistent fields:
\begin{equation}\label{eq:isotropic-magnetic-obstruction}
f(\bx,\bv,0)=F(|\bv|^2),\qquad
\int_{\R^3}F(|\bv|^2)\,d\bv=\ion,\qquad
\bE(\bx,0)=\bB(\bx,0)=0.
\end{equation}
With $\bH=0$, the magnetic transport vanishes because
$\nabla_{\bv}F(|\bv|^2)$ is parallel to $\bv$, and the current
is zero by isotropy. Hence $f(\bx,\bv,t)=F(|\bv|^2)$ and
$\bE=\bB=0$ remain a solution for every admissible external
magnetic field $\bBext(\bx,t)$ in the control model
\eqref{eq:VM-control-3D3V}. This state cannot be driven toward
a different target distribution by purely magnetic actuation.

These observations suggest that magnetic actuation alone is generally insufficient to stabilize the kinetic state. While magnetic control plays an important role in regulating electromagnetic dynamics and particle confinement, an additional mechanism is needed to modify the particle energy distribution and suppress the growth of perturbations. This motivates a combined electromagnetic feedback strategy, in which electric control directly influences the particle energy distribution, while magnetic control regulates the transverse electromagnetic dynamics. Accordingly, throughout this work, we design feedback laws for both the external electric field $\bH$ and the external magnetic field $\bBext$.


\section{Stability analysis and structure-informed control}
In this section, we develop a structure-informed feedback framework for suppressing instabilities in the Vlasov--Maxwell system. We first analyze the linearized longitudinal and transverse dynamics to identify macroscopic quantities relevant to instability growth. We then use the nonlinear perturbation balance to construct a cancellation--dissipation feedback law and examine how delayed observations introduce errors into the field cancellation. The resulting error identities guide the choice of current and field-gradient inputs for a low-rank neural operator that provides a learned correction to the delayed feedback.

For simplicity, we focus on the one-dimensional,
two-velocity (1D--2V) Vlasov--Maxwell system for the controller
construction, subsequent analysis, and numerical experiments. More specifically, $f=f(x,v_1,v_2,t)$, the spatial variable
$x$ lies in a periodic interval, and $\bv=(v_1,v_2)\in\R^2$.
In vector identities involving cross products, we identify
$\bv$ with $(v_1,v_2,0)$ and $\nabla_{\bv}$ with
$(\partial_{v_1},\partial_{v_2},0)$.
The self-consistent fields have the form
$\bE=(E_1,E_2,0)$ and $\bB=(0,0,B_3)$, and the uncontrolled
system is
\begin{equation}\label{eq: VM 1D-2V}
\left\{
\begin{aligned}
&\partial_t f+v_1\partial_x f
 +(E_1+v_2B_3)\partial_{v_1}f
 +(E_2-v_1B_3)\partial_{v_2}f=0,\\
&\partial_xE_1=\rho-\ion,\\
&\partial_tE_1=-J_1,\\
&\partial_tE_2=-J_2-\partial_xB_3,\\
&\partial_tB_3=-\partial_xE_2.
\end{aligned}
\right.
\end{equation}
Here
\begin{equation*}
\rho(x,t)=\int_{\R^2}f(x,\bv,t)\,d\bv,\qquad
J_i(x,t)=\int_{\R^2}v_i f(x,\bv,t)\,d\bv,\quad i=1,2.
\end{equation*}
For zero-mean longitudinal fields, one may also write
$E_1=-\partial_x\Phi$ with a periodic potential satisfying
$-\partial_{xx}\Phi=\rho-\ion$.
The external controls are
$\bH=(H_1,H_2,0)$ and $\bBext=(0,0,\Bext)$.
Their contributions to the Vlasov equation are obtained by
replacing $E_i$ by $E_i+H_i$ and $B_3$ by $B_3+\Bext$ in
\eqref{eq: VM 1D-2V}. The general perturbation equation and
energy balance in Section~2 apply to this reduction with spatial
integrals over the periodic interval and velocity integrals over
$\R^2$.

\subsection{Macroscopic observations for feedback control} \label{sec:linear stability}
Designing a feedback controller based on partial observations requires
identifying macroscopic quantities that convey information about the
coupled kinetic and electromagnetic dynamics. The feedback controller
in \cite{lu2025dynamical} uses the full phase-space perturbation
$\df=f-\feq$ as its input. Here, we seek observations that can be
used to evaluate the feedback without reconstructing the full
distribution. Fourier--Laplace analysis characterizes the amplification
of small perturbations, but the resulting spectral information alone
does not specify which macroscopic quantities should be supplied to
a feedback controller. We first recall the longitudinal response to
clarify this limitation, and then integrate Maxwell's equations to
expose the dependence of the fields on macroscopic currents.

We consider perturbations around a spatially homogeneous equilibrium
$\feq(\bv)$, with $\bEeq=\mathbf 0$, $\Beq=0$, and
$\mathbf J_{eq}=0$. Linearizing the perturbation equation
\eqref{eq: VM perturbation} and retaining first-order contributions
from the control fields yields
\begin{equation}\label{eq: VM linearize}
\left\{
\begin{split}
&\partial_t\df+v_1\partial_x \df+\big(E_1+H_1+v_2(B_3+\Bext)\big)\partial_{v_1}\feq+\big(E_2+H_2-v_1(B_3+\Bext)\big)\partial_{v_2}\feq = 0,\\
& \partial_x E_1 = \drho, \quad \drho = \int \df d\bv,\\
& \partial_t E_1 = -J_1, \quad J_1 = \int v_1\df d\bv, \\
& \partial_t E_2 = -J_2-\partial_x B_3,  \quad J_2 = \int v_2\df d\bv, \\
& \partial_t B_3 = -\partial_x E_2.
\end{split}
\right.
\end{equation}
Following \cite{einkemmer2024control}, we use Fourier and Laplace
transforms to describe the longitudinal linear response. We adopt the
conventions
\begin{equation*}
\widehat{F}(\xi) = \int_{-\infty}^{\infty} F(x)e^{-i\xi x}\,\rd x,
\qquad
\mathcal{L}[F](s) = \int_0^\infty F(t)e^{-st}\,\rd t ,
\end{equation*}
for the Fourier and Laplace transforms, respectively.

For a nonzero spatial Fourier mode, the wavevector is
$(\xi,0,0)$; longitudinal and transverse refer to directions parallel
and perpendicular to this vector, respectively. Longitudinal
instabilities, such as the two-stream and bump-on-tail instabilities,
can appear in the density perturbation $\drho$ and the longitudinal
electric field perturbation $\dE_1$. Applying the Fourier and Laplace
transforms to \eqref{eq: VM linearize} yields
\begin{equation}\label{eq:density perturb}
\begin{split}
\laplace[\widehat{\drho}(\xi,\cdot)](s) = \laplace[\widehat{\drho}_F(\xi,\cdot)](s)&-\frac{\laplace[\widehat{\drho}_F(\xi,\cdot)](s)+i\xi\laplace[\hH_1(\xi,\cdot)](s)}{1+\laplace[U(\xi,\cdot)](s)}\laplace[\widehat{U}(\xi,\cdot)](s)\\
&-i\xi\frac{\laplace[\hB_3(\xi,\cdot)](s)+\laplace[\hB_{ext}(\xi,\cdot)](s)}{1+\laplace[U(\xi,\cdot)](s)}\laplace[\widehat{U_v}(\xi,\cdot)](s).
\end{split}
\end{equation}
Here $\drho_F(x,t) = \int \df_F(x,\bv,t) \rd \bv $
denotes the free-streaming density perturbation, where $\df_F(x,\bv,t) = \df_0(x-v_1 t,v_1,v_2)$ satisfies the longitudinal free-streaming transport equation $\partial_t \df_F+v_1\partial_x \df_F = 0$. The functions $U(x,t)$ and $U_v(x,t)$ are defined through
\[
\widehat{U}(\xi,t) = t\int \widehat{\feq}(t\xi,v_2)dv_2, \qquad \widehat{U_v}(\xi,t) = t\int v_2 \widehat{\feq}(t\xi,v_2) dv_2,
\]
where $\widehat{\feq}$ denotes the Fourier transform with respect to $v_1$. The derivation of \eqref{eq:density perturb} follows standard arguments and is provided in Appendix \ref{append:fourier laplace}.

The longitudinal response in \eqref{eq:density perturb} involves
the function
\[
\mathcal R(\xi,s)=\frac{1}{1+\laplace[\widehat{U}(\xi,\cdot)](s)}.
\]
If there exist $\xi_0\in\R$ and $s_0$ with $\Re(s_0)>0$ such that
\[
1+\laplace[\widehat{U}(\xi_0,\cdot)](s_0)=0,
\]
then, provided that the corresponding numerator in
\eqref{eq:density perturb} does not vanish to sufficient order at
$s_0$, $\laplace[\widehat{\delta\rho}]$ has a pole in the right
half-plane. The inverse Laplace transform then contains a growing
contribution proportional to $e^{s_0t}$. The term involving
$\laplace[\widehat U_v]$ also allows the magnetic response to influence
the density evolution when this coupling kernel is nonzero.

This representation describes how growth can appear in the
longitudinal response and how the control fields enter that response.
It does not, by itself, determine an appropriate set of feedback
observations. In particular, $B_3$ remains a self-consistent field
in \eqref{eq:density perturb}, and the formula does not explicitly
identify the kinetic moment that drives the transverse field dynamics.
Examining the poles of $\mathcal R$ therefore leaves a key part of
the observation design unresolved.

To make this dependence explicit, we now integrate the transverse
Maxwell equations in physical space and time. This provides a direct
representation of the fields in terms of the macroscopic current,
which is more useful here for selecting feedback observations.
The Maxwell equations retain the same form in the full nonlinear
kinetic system, so the following identities are not restricted to
the linearized dynamics. Introducing the characteristic variables
\[
W^+ = E_2+B_3, \qquad W^- = E_2-B_3,
\]
we rewrite the last two equations of \eqref{eq: VM linearize} as
\[
\begin{split}
\partial_t W^+ +\partial_x W^+ = -J_2, \quad \partial_t W^- - \partial_x W^- = -J_2.
\end{split}
\]
Integrating these equations along the characteristic curves $\dot{x} = 1$ and $\dot{x}=-1$, respectively, gives
\[
\begin{split}
W^+(x,t) &= W^+_0(x-t)-\int^t_0 J_2(x-t+s,s)d s,\\
W^-(x,t) &= W^-_0(x+t)-\int^t_0 J_2(x+t-s,s)d s.
\end{split}
\]
The transverse electric and magnetic fields can then be recovered from $\ds E_2 = \frac{W^+ +W^-}{2}$ and $\ds B_3 = \frac{W^+ -W^-}{2}$, yielding
\begin{equation}\label{eq:E2 B3}
\begin{split}
E_2(x,t) &= \frac{1}{2}\Big( E_2(x+t,0) + E_2(x-t,0) \Big)-\frac{1}{2}\Big(B_3(x+t,0)-B_3(x-t,0)\Big)
\\
&\quad -\frac{1}{2}\int^{t}_{0}J_2(x-t+s,s)\,d s
      -\frac{1}{2}\int^{t}_{0}J_2(x+t-s,s)\,d s,\\
B_3(x,t) & = \frac{1}{2}\Big( B_3(x+t,0) + B_3(x-t,0) \Big) - \frac{1}{2}\Big(E_2(x+t,0)-E_2(x-t,0)\Big)
\\
&\quad + \frac{1}{2}\int_0^t  J_2(x+t-s,s) d s-\frac{1}{2}\int^t_0 J_2(x-t+s,s)  d s .
\end{split}
\end{equation}
Equation \eqref{eq:E2 B3} explicitly separates the transverse field
evolution into freely propagating initial fields and contributions
from the history of $J_2$. It identifies the transverse current as
the moment through which the kinetic distribution enters the
transverse Maxwell dynamics. In particular, a growing current mode
can generate corresponding growth in $E_2$ and $B_3$. The value of
this representation for feedback design is that it exposes the
current--field dependence directly; it does not require a separate
transverse dispersion calculation to identify $J_2$ as a relevant
observation.

The integral representation thus supplies information about the
field dynamics that the longitudinal pole criterion alone does not
provide. Together with the density--electric-field relation in the
longitudinal direction, it motivates the macroscopic quantities
$\rho$, $E_1$, $J_2$, $E_2$, and $B_3$ as candidate feedback
observations. These relations do not yet prescribe how the control
fields should depend on the observations or establish that the
observations suffice for stabilization. In the next subsection,
we use the nonlinear perturbation balance to construct a feedback
law and then analyze the cancellation errors caused by delayed
observations. The resulting identities further identify the currents
and spatial field derivatives used as inputs to the learned
correction.

\subsection{A cancellation-dissipation structure and time-delay compensation}
Equations \eqref{eq:density perturb}--\eqref{eq:E2 B3} reveal how the self-consistent electromagnetic fields are coupled to moments of the distribution function perturbation $\delta f$. These observations suggest that feedback control based on the macroscopic quantities $\rho$, $E_1$, $J_2$, $E_2$, and $B_3$ may be effective in suppressing both longitudinal and transverse instabilities. However, the appropriate functional dependence of the control fields on these quantities remains unclear. Moreover, the complex coupling between kinetic and electromagnetic dynamics makes it challenging to derive a theoretically optimal feedback control law.

For practical implementation, we instead pursue a suboptimal feedback control strategy with a simple structure. Motivated by the perturbation energy identity \eqref{eq:energy estimate}, derived from the \emph{original nonlinear} Vlasov--Maxwell system, we consider a linear feedback control law of the form
\begin{equation}\label{eq:cancellation based}
\begin{split}
H_1[\df] &= -\dE_1+\gam_1\int\df\partial_{v_1}\feq d\bv,\\
H_2[\df] &= -\dE_2+\gam_2\int\df\partial_{v_2}\feq d\bv,\\
\Bext[\df] &= -\dB_3-\gam_b\int(v_1\partial_{v_2}\feq-v_2\partial_{v_1}\feq)\df d\bv\,,
\end{split}
\end{equation}
which yields the following dissipation identity:
\begin{align*}
& \frac{1}{2}\frac{\rd}{\rd t }||\df(t)||^2_2
\\ &=-\gam_1\int|\int\df\partial_{v_1}\feq d\bv|^2 \rd x-\gam_2\int|\int\df\partial_{v_2}\feq d\bv|^2 \rd x-\gam_b\int|\int(v_1\partial_{v_2}\feq-v_2\partial_{v_1}\feq)\df d\bv|^2 \rd x\leq 0\,.
\end{align*}
In particular, for sufficiently small initial perturbations, pure field cancellation ($\gam_1 = \gam_2 = \gam_b = 0$),
\[
H_1 = -\dE_1, \quad H_2 = -\dE_2, \quad \Bext = -\dB_3,
\]
provides a neutral stabilization mechanism, satisfying $\frac{\rd}{\rd t }||\df(t)||^2_2 = 0$. In practice, however, the controller may not respond instantaneously to changes in the system due to factors such as the computational time required by the numerical solver and delays in adjusting the control fields. Such delayed responses can significantly degrade the effectiveness of pure cancellation control, as demonstrated in Section \ref{sec:experiments}.

To account for the effects of time delay, we assume that the control applied at time $t$ is computed using information available at an earlier time $t-\tau$ (or before), where $\tau>0$ denotes the time lag. The resulting cancellation errors can be expressed as
\begin{equation}\label{eq:cancellation err}
\begin{split}
\dE_1(x,t)-\dE_1(x,t-\tau) &= \int^t_{t-\tau}\partial_t\dE_1(x,s)ds = -\int^t_{t-\tau}\dJ_1(x,s)ds \\
\dE_2(x,t)-\dE_2(x,t-\tau) &= \int^t_{t-\tau}\partial_t\dE_2(x,s)ds= -\int^t_{t-\tau}\big(\dJ_2(x,s)+\partial_x\dB_3(x,s)\big) ds\\
\dB_3(x,t)-\dB_3(x,t-\tau) &= \int^t_{t-\tau}\partial_t\dB_3(x,s)ds = -\int^t_{t-\tau}\partial_x\dE_2(x,s)ds
\end{split}
\end{equation}
Motivated by these identities, we introduce operator-based dissipative corrections, $\dH_1$, $\dH_2$, and $\delta \Bext$, that utilize the observations $\dJ_1(t-\tau)$, $\dJ_2(t-\tau)$, $\partial_x\dE_2(t-\tau)$, and $\partial_x\dB_3(t-\tau)$ to compensate for delay-induced cancellation errors. Specifically, we propose the following cancellation-residual feedback control structure:
\begin{equation}\label{eq:feedback operator}
\begin{split}
H_1(x,t) &= -\delta E_1(x,t-\tau)+\network_{\th_1}[\dJ_1(t-\tau),\dJ_2(t-\tau),\partial_x \dE_2(t-\tau),\partial_x \dB_3(t-\tau)](x),\\
H_2(x,t) &= -\delta E_2(x,t-\tau)+\network_{\th_2}[\dJ_1(t-\tau),\dJ_2(t-\tau),\partial_x \dE_2(t-\tau), \partial_x\dB_3(t-\tau)](x),\\
\Bext(x,t) &= -\delta B_3(x,t-\tau)+\network_{\th_3}[\dJ_1(t-\tau),\dJ_2(t-\tau), \partial_x \dE_2(t-\tau), \partial_x \dB_3 (t-\tau)](x).
\end{split}
\end{equation}
Here, $\network_{\th_i}$ are neural operators that map the observed quantities to the corresponding control corrections. Possible architectures include the Fourier Neural Operator (FNO) \cite{li2020fourier,kovachki2021universal,li2023fourier}, Low-Rank Neural Operator (LNO) \cite{kovachki2023neural}, Graph Neural Operator (GNO) \cite{anandkumar2020neural,li2020multipole}, and DeepONet \cite{lu2021learning,wang2021learning,lanthaler2022error}. Although selecting an optimal operator architecture is nontrivial, motivated by the linear dissipation structure of the feedback law \eqref{eq:cancellation based}, we propose a linear neural operator with a low-rank structure:
\begin{equation}\label{eq:dissipation}
\begin{split}
&\network_\th [\dJ_1(t),\dJ_2(t), \partial_x \dE_2(t), \partial_x \dB_3(t)](x) \\&= \sum^r_{k=1} \phi_k(x)\big[\int\psi^{(1)}_k(y;\th)\dJ_1(y,t) dy+\int\psi^{(2)}_k(y;\th)\dJ_2(y,t) dy\\
&\quad +\int\psi^{(3)}_k(y;\th)\partial_x \dE_2(y,t) dy+\int\psi^{(4)}_k(y;\th)\partial_x \dB_3(y,t) dy\big].
\end{split}
\end{equation}
Here, $\{\phi_k(x)\}_{k\geq 0}$ forms an $L^2$ basis on the computational domain, which can either be learned or chosen analytically (e.g., a trigonometric basis under periodic boundary conditions). The kernels $\{\psi^{(j)}_k(x;\th)\}_{1\leq j\leq 4, 1\leq k \leq r}$ are parameterized by neural networks. The structure in \eqref{eq:dissipation} resembles a DeepONet-type architecture with a linear branch network, where the branch coefficients are computed through kernel integrals rather than a conventional branch network, allowing greater flexibility in handling different computational meshes.

In summary, the proposed feedback controller \eqref{eq:feedback operator} integrates three key components:
\begin{itemize}
    \item Instability-informed macroscopic observations identified through the linear stability analysis;
    \item A field-cancellation baseline motivated by the energy structure of the Vlasov--Maxwell system;
    \item Delay-informed current and field-gradient observations, together with learned residual corrections motivated by Maxwell's equations.
\end{itemize}

\section{Training of the control network}\label{sec:training}
In this section, we describe the training procedure for the neural feedback controller. The network parameters are optimized using the adjoint-state method within an optimize-then-discretize framework. We also introduce a semi-Lagrangian discretization of the kinetic equations to enable efficient long-time simulations.

\subsection{Adjoint state method for controller training}

Training the feedback operator \eqref{eq:feedback operator} requires solving the following PDE-constrained optimization problem over the time interval $[0,T]$:
\begin{subequations}\label{eq:fb optim}
\begin{equation}\label{eq:fb optim loss}
\min_{\th = (\th_1,\th_2,\th_3)} \quad \cJ(\th) = \frac{1}{2}\int^T_0||f(\cdot,\cdot,t;H_1(\th_1),H_2(\th_2),\Bext(\th_3))-\feq(\cdot,\cdot)||_{x,\bv}^2 d t
\end{equation}
\begin{equation}\label{eq:fb optim constraint}
s.t. \quad  \left\{
\begin{array}{l}
\partial_t f +v_1\partial_x f +(\bE+\bH)\cdot\nabla_{\bv}f+\bv\times(\bB+\bBext)\cdot\nabla_{\bv}f = 0\,, \\
\\
E_1 = - \partial_x \Phi, \quad \partial_x E_1 = -\partial_{xx}\Phi = \rho-\ion, \\
\\
\partial_t E_1 = -J_1, \quad J_1 = \int v_1 f d\bv, \\
\\
\partial_t E_2 = -J_2-\partial_x B_3, \quad J_2 = \int v_2 f d\bv \\
\\
\partial_t B_3 = -\partial_x E_2,
\end{array}
\right.
\end{equation}
\end{subequations}
To obtain an effective control operator, the optimization time horizon must be sufficiently long to capture the growth of instabilities while remaining computationally manageable. Moreover, the controller should be exposed to diverse training trajectories to promote its effectiveness beyond the training interval $[0,T]$. To this end, we incorporate a variety of initial perturbations during training, as described in Section \ref{sec:experiments}.

We solve the optimization problem \eqref{eq:fb optim} using the adjoint-state method. To this end, we introduce the Lagrangian
\[
\cL(\lam,\th) = \frac{1}{2}\int^T_0||\df(t)||^2_2\rd t +\int^T_0\iint\lambda(x,\bv,t)\big(\partial_t f +v_1\partial_x f +(\bE+\bH)\cdot\nabla_{\bv}f+\bv\times(\bB+\bBext)\cdot\nabla_{\bv}f\big)\rd xd\bv \rd t \,.
\]
When $f$ satisfies the constraint \eqref{eq:fb optim constraint}, the constraint term vanishes, so $\cL(\lam,\th) = \cJ(\th)$ for arbitrary $\lam$ and $\th$, and consequently $\nabla_\th \cJ = \nabla_\th \cL$. The resulting gradient can then be used in standard optimization algorithms, such as SGD or Adam, to update the control parameters.

Applying integration by parts yields
\[
\begin{split}
\nabla_\th\cJ = \nabla_\th\cL &= \int^T_0\iint \df(x,\bv,t)\nabla_\th f \rd xd\bv \rd t +\iint\lam(x,\bv,t)\nabla_\th f(x,\bv,t)|^T_0 \rd x  d\bv\\
&-\int^T_0\iint(\partial_t\lam+v_1\partial_x\lam+(\bE+\bH)\cdot\nabla_\bv\lam+\bv\times(\bB+\bBext)\cdot\nabla_\bv\lam)\nabla_\th f \rd x  d\bv \rd t \\
&+\int^T_0\iint\big(\nabla_\th\bE(x,t)+\nabla_\th(\bH(x,t))\big)\cdot\nabla_\bv f(x,\bv,t) \lam(x,\bv,t) \rd x \rd \bv  \rd t \\
&-\int^T_0\iint\big(\nabla_\th\bB+\nabla_\th(\bBext(x,t))\big)\cdot\big(\bv\times\nabla_\bv f(x,\bv,t)\big)\lam(x,\bv,t)\rd x d\bv \rd t 
\end{split}
\]
If $\lambda$ satisfies the adjoint equation
\begin{equation}\label{eq:adjoint pde}
\partial_t\lam+v_1\partial_x\lam+(\bE+\bH)\cdot\nabla_\bv\lam+\bv\times(\bB+\bBext)\cdot\nabla_\bv\lam = \df, \quad \lambda(x,v,T) = 0,
\end{equation}
the gradient of the objective function simplifies to
\begin{equation}\label{eq:grad}
\begin{split}
\nabla_\th\cJ &= \int^T_0\iint\big(\nabla_\th\bE(x,t)+\nabla_\th(\bH(x,t))\big)\cdot\nabla_\bv f(x,\bv,t) \lam(x,\bv,t) \rd x \rd \bv  \rd t \\
&-\int^T_0\iint\big(\nabla_\th\bB+\nabla_\th(\bBext(x,t))\big)\cdot\big(\bv\times\nabla_\bv f(x,\bv,t)\big)\lam(x,\bv,t)\rd x d\bv \rd t .
\end{split}
\end{equation}
Substituting the feedback operator \eqref{eq:feedback operator} into the above expression, we obtain
\[
\begin{split}
\nabla_{\th}(H_1(x,t)) =& -\nabla_\th E_1(x,t-\tau)+\network_{\th_1}[\nabla_\th J_1(t-\tau),\nabla_\th J_2(t-\tau),\partial_x \nabla_\th E_2(t-\tau),\partial_x \nabla_\th B_3(t-\tau)](x)+\\
&(\nabla_{\th_1}\network_{\th_1})[\dJ_1(t-\tau),\dJ_2(t-\tau),\partial_x \dE_2(t-\tau),\partial_x \dB_3(t-\tau)](x),
\end{split}
\]

\[
\begin{split}
\nabla_{\th}(H_2(x,t)) =& -\nabla_\th E_2(x,t-\tau)+\network_{\th_2}[\nabla_\th J_1(t-\tau),\nabla_\th J_2(t-\tau),\partial_x \nabla_\th E_2(t-\tau),\partial_x \nabla_\th B_3(t-\tau)](x)+\\
&(\nabla_{\th_2}\network_{\th_2})[\dJ_1(t-\tau),\dJ_2(t-\tau),\partial_x \dE_2(t-\tau),\partial_x \dB_3(t-\tau)](x),
\end{split}
\]

\[
\begin{split}
\nabla_{\th}(\Bext(x,t)) =& -\nabla_\th B_3(x,t-\tau)+\network_{\th_3}[\nabla_\th J_1(t-\tau),\nabla_\th J_2(t-\tau),\partial_x \nabla_\th E_2(t-\tau),\partial_x \nabla_\th B_3(t-\tau)](x)+\\
&(\nabla_{\th_3} \network_{\th_3})[\dJ_1(t-\tau),\dJ_2(t-\tau),\partial_x \dE_2(t-\tau),\partial_x \dB_3(t-\tau)](x),
\end{split}
\]
where, for $1\leq i\leq 3$, the operator terms can be written explicitly as
\[
\begin{split}
\network_{\th_i}[\nabla_\th J_1,\nabla_\th J_2,\nabla_\th E_2, \nabla_\th B_3](x) & = \sum^r_{k=1} \phi_k(x)\big[\int \psi^{(1)}_k(y;\th_i)\nabla_\th J_1(y) \rd y +\int\psi^{(2)}_k(y;\th_i)\nabla_\th J_2(y) \rd y \\
&\quad +\int\psi^{(3)}_k(y;\th_i)\partial_x \nabla_\th E_2(y) \rd y +\int\psi^{(4)}_k(y;\th_i)\partial_x \nabla_\th B_3(y) \rd y\big].
\end{split}
\]
and
\[
\begin{split}
(\nabla_{\th_i} \network_{\th_i})[\dJ_1,\dJ_2,\partial_x \dE_2,\partial_x \dB_3](x)& = \sum^r_{k=1} \phi_k(x)\big[\int\nabla_{\th_i}\psi^{(1)}_k(y;\th_i)\dJ_1(y) dy+\int\nabla_{\th_i}\psi^{(2)}_k(y;\th_i)\dJ_2(y) dy\\
&\quad +\int\nabla_{\th_i}\psi^{(3)}_k(y;\th_i)\partial_x \dE_2(y) dy+\int\nabla_{\th_i}\psi^{(4)}_k(y;\th_i)\partial_x \dB_3(y) dy\big].
\end{split}
\]

In general, the state sensitivities $\nabla_\th J_1$, $\nabla_\th J_2$, $\nabla_\th E_1$, $\nabla_\th E_2$, and $\nabla_\th B_3$ are not available in closed form, and computing them numerically through automatic differentiation can be expensive. To reduce the computational cost, we treat the state-dependent observations as fixed during each gradient evaluation, leading to the approximate gradient
\begin{equation}\label{eq:grad approx}
\begin{split}
\nabla_\th\cJ\approx&\int^T_0 \iint (\nabla_{\th_1}\network_{\th_1})[\dJ_1(t-\tau),\dJ_2(t-\tau),\partial_x\dE_2(t-\tau),\partial_x\dB_3(t-\tau)] \cdot\lam \partial_{v_1} f \rd x \rd\bv \rd t\\
&+\int^T_0 \iint (\nabla_{\th_2}\network_{\th_2})[\dJ_1(t-\tau),\dJ_2(t-\tau),\partial_x\dE_2(t-\tau),\partial_x\dB_3(t-\tau)]\cdot \lam\partial_{v_2}f \rd x \rd\bv \rd t\\
&+\int^T_0 \iint (\nabla_{\th_3}\network_{\th_3})[\dJ_1(t-\tau),\dJ_2(t-\tau),\partial_x\dE_2(t-\tau),\partial_x\dB_3(t-\tau)]\cdot\lam(v_2\partial_{v_1}f-v_1\partial_{v_2}f) \rd x \rd\bv \rd t
\end{split}
\end{equation}

We emphasize that the approximation \eqref{eq:grad approx} leads to a sequential optimization procedure: the state trajectory is held fixed during each gradient evaluation and recomputed after each parameter update. The exact closed-loop gradient \eqref{eq:grad} contains two contributions. The first arises from the explicit dependence of the learned control on its parameters along the current trajectory, while the second accounts for changes in the feedback observations induced by variations in the state trajectory. The approximation \eqref{eq:grad approx} retains only the first contribution during each gradient evaluation. Rather than explicitly differentiating the second contribution, its effects are incorporated indirectly through the updated trajectories in subsequent iterations.

This approach avoids differentiating through the entire closed-loop dynamics while retaining the direct sensitivity of the objective function to the learned control. It is particularly suitable for the proposed cancellation-residual formulation, in which the learned dissipation operator provides a moderate correction to a prescribed field-cancellation baseline. Although the approximation sacrifices the exact closed-loop gradient in favor of computational efficiency, our numerical results indicate that the sequential optimization procedure produces stable parameter updates when the initial perturbation is sufficiently small and an appropriate optimization step size is used. Its performance will be further evaluated through numerical experiments in Section~\ref{sec:experiments}.

\subsection{Semi-Lagrangian solver}

Solving the PDE-constrained optimization problem \eqref{eq:fb optim} requires repeatedly solving high-dimensional kinetic equations. To relax the time-step restrictions and improve the efficiency of long-time simulations, we employ a semi-Lagrangian scheme for the kinetic equation, coupled with a leapfrog-type update for the electromagnetic fields. Let $f^n$, $E^n_1$, $E^n_2$, and $B^n_3$ denote the numerical solutions at time $t=t^n$. We assume that $f^n$ is represented on the space-velocity mesh $\{(x_i,v_{1,j},v_{2,k})\}$, $E^n_1$ and $E^n_2$ on the spatial mesh $\{x_i\}$, and $B^n_3$ on the staggered mesh $\{x_{\iph}\}$, where $\ds x_{\iph} = \frac{1}{2}(x_i+x_{i+1})$. The numerical solution is advanced from $t_n$ to $t_{n+1}$ through the following steps:

\begin{itemize}
    \item First, compute $f^{(1)}(x_i,v_{1,j},v_{2,k}) = f^n(x_i-\frac{\Dt}{2}v_{1,j},v_{1,j},v_{2,k})$.

    \item Compute $E_1^{\nph}$ from $\rho^{(1)} = \int f^{(1)}d\bv$ and $J_2^{\nph} = \int v_2 f^{(1)}d\bv$. Then update $E_2$ on the primal mesh $\{x_i\}$ and $B_3$ on the staggered mesh $\{x_{\iph}\}$ according to
    \[
    \begin{split}
    &\frac{B^{\nph}_{3,\iph}-B^n_{3,\iph}}{(\Dt/2)} = -\frac{E^n_{2,i+1}-E^n_{2,i}}{\Dx},\\
    &\frac{E^{\nph}_{2,i}-E^n_{2,i}}{(\Dt/2)} = -J^{\nph}_{2,i}-\frac{B^{\nph}_{3,\iph}-B^{\nph}_{3,\imh}}{\Dx}.\\
    \end{split}
    \]
    Next, obtain the nodal values of $B_3$ on the primal mesh by averaging, $\ds \baB^{\nph}_{3,i} = \frac{1}{2}(B^{\nph}_{3,\imh}+B^{\nph}_{3,\iph})$, and use these values to compute
    \[
    \begin{split}
    f^{(2)}(x_i,v_{1,j},v_{2,k}) = f^{(1)}\big(x_i,&v_{1,j}-\Dt[E_{1,i}^{\nph}+H^{\nph}_{1,i}+v_{2,k}(\baB^{\nph}_{3,i}+B^{\nph}_{ext,i})], \\
    &v_{2,k}-\Dt[E^{\nph}_{2,i}+H^{\nph}_{2,i}-v_{1,j}(\baB^{\nph}_{3,i}+B^{\nph}_{ext,i})]\big)
    \end{split}
    \]

    \item Finally, compute $f^{n+1}(x_i,v_{1,j},v_{2,k}) = f^{(2)}(x_i-\frac{\Dt}{2}v_{1,j},v_{1,j},v_{2,k})$. Obtain $E_1^{n+1}$ from $\rho^{n+1} = \int f^{n+1}d\bv$ and complete the updates of $E_2$ and $B_3$ using
    \[
    \begin{split}
    &\frac{E^{n+1}_{2,i}-E^{\nph}_{2,i}}{(\Dt/2)} = -J^{\nph}_{2,i}-\frac{B^{\nph}_{3,\iph}-B^{\nph}_{3,\imh}}{\Dx},\\
    &\frac{B^{n+1}_{3,\iph}-B^{\nph}_{3,\iph}}{(\Dt/2)} = -\frac{E^{n+1}_{2,i+1}-E^{n+1}_{2,i}}{\Dx}.
    \end{split}
    \]
\end{itemize}

The updates of $E_2$ and $B_3$ can be expressed in the following leapfrog form:
\begin{equation}\label{eq:yee's lattice}
\begin{split}
&\frac{B^{\nph}_{3,\iph}-B^{\nmh}_{3,\iph}}{\Dt} = -\frac{E^n_{2,i+1}-E^n_{2,i}}{\Dx}, \\
&\frac{E^{n+1}_{2,i}-E^n_{2,i}}{\Dt} = -J^{\nph}_{2,i}-\frac{B^{\nph}_{3,\iph}-B^{\nph}_{3,\imh}}{\Dx},
\end{split}
\end{equation}
which recovers the classical Yee scheme \cite{yee1966numerical}. The scheme \eqref{eq:yee's lattice} provides a second-order time discretization of the semi-discrete system
\begin{equation}\label{eq:yee's semidiscrete}
\begin{split}
&\frac{\rd}{\rd t}B_{3,\iph}(t) = -\frac{1}{\Dx}(E_{2,i+1}(t)-E_{2,i}(t)), \\
&\frac{\rd}{dt}E_{2,i}(t) = -J_{2,i}(t)-\frac{1}{\Dx}(B_{3,\iph}(t)-B_{3,\imh}(t)).
\end{split}
\end{equation}
An important property of the semi-discrete scheme \eqref{eq:yee's semidiscrete} is that it satisfies a discrete energy balance. To see this, define the transverse electromagnetic energy by
\[
U(t) = \frac{1}{2}\sum_i (|E_{2,i}(t)|^2+|B_{3,\iph}(t)|^2) \Dx
\]
Differentiating with respect to time and substituting \eqref{eq:yee's semidiscrete}, we obtain
\begin{equation}\label{eq:yee's energy est1}
\begin{split}
\frac{\rd}{dt}U(t) &= \sum_i E_{2,i}\frac{\rd}{dt}E_{2,i}(t)\Dx+\sum_i B_{3,\iph}\frac{\rd}{dt}B_{3,\iph}(t)\Dx\\
& = -\sum_i J_{2,i}E_{2,i}\Dx-\sum_i E_{2,i}(B_{3,\iph}-B_{3,\imh})-\sum_i B_{3,\iph}(E_{2,i+1}-E_{2,i}).
\end{split}
\end{equation}
Using the discrete summation-by-parts identity
\[
\sum_i E_{2,i}(B_{3,\iph}-B_{3,\imh}) = -\sum_i (E_{2,i+1}-E_{2,i})B_{3,\iph},
\]
we see that the last two terms on the right-hand side of \eqref{eq:yee's energy est1} cancel. Consequently,
\[
\frac{\rd}{dt}U(t) = -\sum_i J_{2,i}E_{2,i}\Dx.
\]
In the vacuum case, where $J_2 = 0$, the energy $U(t)$ is conserved. However, this energy conservation property does not generally carry over to the fully discrete scheme \eqref{eq:yee's lattice}. To achieve energy conservation at the fully discrete level, implicit time integrators, such as the Crank--Nicolson scheme, may be employed.

The adjoint equation \eqref{eq:adjoint pde} can be integrated backward in time using a similar semi-Lagrangian procedure:

\begin{itemize}
    \item First, compute $\lambda^{(-1)}(x_i,v_{1,j},v_{2,k}) = \lambda^{n}(x_i+\frac{\Dt}{2}v_{1,j},v_{1,j},v_{2,k})$.
    \item Next, compute
    \[
    \begin{split}
    \lambda^{(-2)}(x_i,v_{1,j},v_{2,k}) = \lambda^{(-1)}\big(x_i,&v_{1,j}+\Dt[E_{1,i}^{\nph}+H^{\nph}_{1,i}+v_{2,k}(\baB^{\nph}_{3,i}+B^{\nph}_{ext,i})], \\
    &v_{2,k}+\Dt[E^{\nph}_{2,i}+H^{\nph}_{2,i}-v_{1,j}(\baB^{\nph}_{3,i}+B^{\nph}_{ext,i})]\big)
    \end{split}
    \]
    \item Finally, compute $\lambda^{n-1}(x_i,v_{1,j},v_{2,k}) = \lambda^{(-2)}(x_i+\frac{\Dt}{2}v_{1,j},v_{1,j},v_{2,k})$.
\end{itemize}

The overall PDE-constrained optimization procedure is summarized in Algorithm \ref{alg:vm-control}.

\begin{algorithm}[h!]
\caption{PDE-constrained optimization for VM control}
\label{alg:vm-control}

Specify the terminal time $T$ and the maximum number of optimization iterations $\texttt{maxiter}$.

Initialize the control parameters
$\theta_1^{0},\theta_2^{0},\theta_3^{0}$.

\begin{algorithmic}[1]

\For{$n=0,\ldots,\texttt{maxiter}-1$}

    \State Solve the controlled VM system \eqref{eq:fb optim constraint} forward in time over $[0,T]$ using the current control laws
    $H_1(\theta_1^{n})$, $H_2(\theta_2^{n})$, and
    $B_{\mathrm{ext}}(\theta_3^{n})$

    \State Solve the adjoint equation
    \eqref{eq:adjoint pde} backward in time to obtain $\lambda^{n}$

    \State Evaluate the approximate gradient using
    \eqref{eq:grad approx}

    \State Update the control parameters
    $\theta_1^{n+1}$, $\theta_2^{n+1}$, and $\theta_3^{n+1}$

\EndFor
\end{algorithmic}
\end{algorithm}

\section{Numerical experiments}\label{sec:experiments}

We evaluate the performance of the electromagnetic feedback control
\eqref{eq:feedback operator} through numerical experiments. In all test cases,
periodic boundary conditions are imposed at both endpoints of the spatial
interval. To construct the dissipation operator \eqref{eq:dissipation}, the
integral kernels
\[
\psi^{NN}(x;\th)
=
\left\{\psi_k^{(j)}(x;\th)\right\}_{1\leq j\leq 4,\;1\leq k\leq r}
\]
are represented by a fully connected multi-output neural network with
architecture $1$-$32$-$32$-$4r$. The output of each hidden layer is activated
by the \texttt{ReLU} function before being passed to the next layer. The
controllers are optimized with respect to the running loss
\[
\frac{1}{2}\int_0^T \|f(t)-\feq\|_2^2\,\rd t,
\qquad T=20,
\]
following the procedure described in Section~\ref{sec:training}.

To stabilize the optimization during the early iterations, we initialize the
parameters in the output layer of $\psi^{NN}$ to zero, so that the learned
dissipation initially vanishes. The network parameters are then updated using
the approximate gradient \eqref{eq:grad approx} for 1000 iterations. We first
perform 200 \texttt{Adagrad} iterations as a preprocessing stage, followed by
\texttt{Adam} iterations for further acceleration.

To expose the controller to a richer family of off-equilibrium states and
improve generalization beyond the optimization horizon $t\in[0,T]$, we add an
artificial perturbation $f_p$ to the initial data and compute the gradient
along the resulting perturbed trajectory at each iteration. The perturbation
takes the form
\[
f_p(x,v_1,v_2)
=
\eps_p
\sum_{k=1}^{K}
\sum_{n_1=1}^{N_1}
\sum_{n_2=1}^{N_2}
\omega_{kn_1n_2}
\hphi_k(x)\hh_{n_1}(v_1)\hh_{n_2}(v_2).
\]
Here $\{\hphi_k(x)\}$, $\{\hh_{n_1}(v_1)\}$, and
$\{\hh_{n_2}(v_2)\}$ are orthonormal bases in the $x$, $v_1$, and $v_2$
directions, respectively. In our computations, we use the normalized
trigonometric functions on $x\in[0,10\pi]$,
\[
\{\hphi_k(x)\}_{k=1}^{11}
=
\left\{
\sqrt{\frac{1}{10\pi}},
\sqrt{\frac{1}{5\pi}}\sin\left(\frac{l}{5}x\right),
\sqrt{\frac{1}{5\pi}}\cos\left(\frac{l}{5}x\right)
\right\}_{1\leq l\leq 5},
\]
and the orthonormal Hermite basis functions on
$v_i\in(-\infty,\infty)$,
\[
\hh_{n_i}(v_i)
=
\sqrt{\frac{1}{\sqrt{2\pi}n_i!}}\,
He_{n_i}(v_i)\exp\left(-\frac{v_i^2}{4}\right),
\qquad
0\leq n_i\leq 5,\quad i=1,2.
\]
The Hermite polynomials are generated recursively by
\[
He_0(v)=1,
\qquad
He_{n+1}(v)=vHe_n(v)-nHe_{n-1}(v),
\qquad n\geq 1.
\]
We set the perturbation magnitude to $\eps_p=10^{-3}$.

The full tensor-product family contains $11\times 6\times 6=396$ basis
functions. To reduce the sampling dimension, we randomly select $M=200$
basis functions from this family and sample their coefficient vector
$\omega\in\R^M$ uniformly from the $M$-dimensional unit ball.

During training, we use the future loss $\frac{1}{2}\|f(2T)-f_{\mathrm{eq}}\|_2^2$
as a measure of long-horizon generalization and select the best-performing
iteration accordingly.

\subsection{Two-stream instability}

We first consider the two-stream equilibrium
\[
\begin{split}
&\feq(v_1,v_2)
=
\frac{1}{4\pi}
\exp\left(
-\frac{(v_1-\bav_1)^2}{2}
-\frac{v_2^2}{2}
\right)
+
\frac{1}{4\pi}
\exp\left(
-\frac{(v_1+\bav_1)^2}{2}
-\frac{v_2^2}{2}
\right),\\
&E_{eq,1}=0,\qquad E_{eq,2}=0,\qquad B_{eq,3}=0,
\end{split}
\]
with $\bav_1=2.4$. The system is initialized with
\[
\begin{split}
&f_0(x,\bv)
=
\left(1+0.005\cos\left(\frac{x}{5}\right)\right)\feq(\bv),
\qquad x\in[0,10\pi],\\
&E_2(x,0)=0,\qquad
B_3(x,0)=0.001\sin\left(\frac{x}{5}\right).
\end{split}
\]
The two-stream equilibrium is unstable, so a small initial perturbation is
amplified as the system evolves. Our goal is to stabilize the full kinetic
distribution around the prescribed equilibrium, with particular attention
to the effect of feedback delay.

\begin{figure}[h!]
    \centering
    \begin{subfigure}{0.4\textwidth}
    \includegraphics[width=\textwidth]{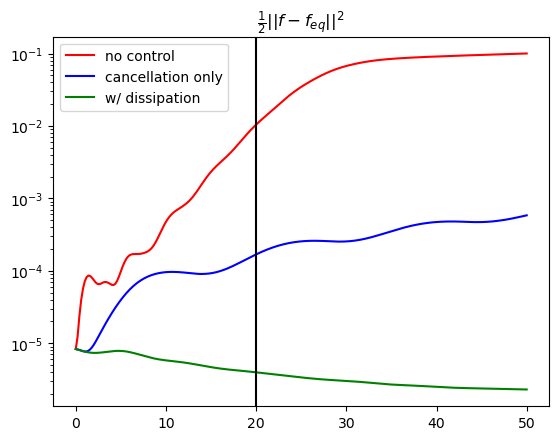}
    \subcaption{$L^2$ perturbation}
    \end{subfigure}
    \quad
    \begin{subfigure}{0.4\textwidth}
    \includegraphics[width=\textwidth]{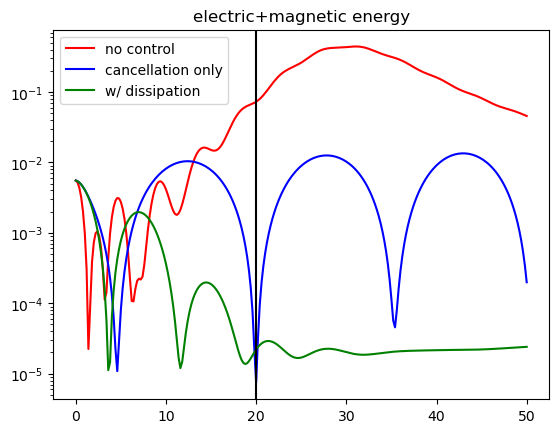}
    \subcaption{Electromagnetic energy}
    \end{subfigure}
    \caption{Two-stream instability. Histories of the $L^2$ perturbation and electromagnetic energy over $t\in[0,50]$. No external field (red) vs.\ delayed cancellation control (blue) vs.\ cancellation--dissipation control (green). The feedback delay is $\tau=\tau_{\mathrm{ref}}=0.6$. The vertical black line marks the terminal time of optimization.}
    \label{fig:two stream time lag perturb}
\end{figure}

We take the reference feedback delay to be $\tau_{\mathrm{ref}}=0.6$ and
optimize the controller over the time interval $t\in[0,20]$. The
computational domain is
\[
(x,v_1,v_2)\in[0,10\pi]\times[-8,8]\times[-8,8].
\]
The learning rate is set to $10^{-4}$ during the \texttt{Adagrad}
preprocessing stage and to $5\times10^{-5}$ during the subsequent
\texttt{Adam} iterations. Figure~\ref{fig:two stream time lag perturb}
compares the uncontrolled solution, delayed cancellation control,
\[
H_1(t)=-\dE_1(t-\tau_{\mathrm{ref}}),\qquad
H_2(t)=-\dE_2(t-\tau_{\mathrm{ref}}),\qquad
\Bext(t)=-\dB_3(t-\tau_{\mathrm{ref}}),
\]
and the cancellation--dissipation feedback control
\eqref{eq:feedback operator}. Without control, both the kinetic
perturbation
\[
\frac{1}{2}\|\df(t)\|_{L^2_{x,v}}^2
\]
and the electromagnetic energy
\[
\energy(t)
=
\frac{1}{2}
\left(
\|E_1(t)\|_{L^2_x}^2
+\|E_2(t)\|_{L^2_x}^2
+\|B_3(t)\|_{L^2_x}^2
\right)
\]
grow rapidly. For instantaneous observations, the pure cancellation law
$H_1=-\dE_1$, $H_2=-\dE_2$, and $\Bext=-\dB_3$ keeps the
$L^2$ perturbation norm constant, as predicted by
\eqref{eq:energy estimate}. Its effectiveness, however, deteriorates when
the cancellation is evaluated using delayed observations. In comparison,
the learned cancellation--dissipation controller substantially improves
the suppression of the instability and remains effective well beyond the
optimization horizon.

\begin{figure}[h!]
    \centering
    \begin{subfigure}{0.3\textwidth}
    \includegraphics[width=\textwidth]{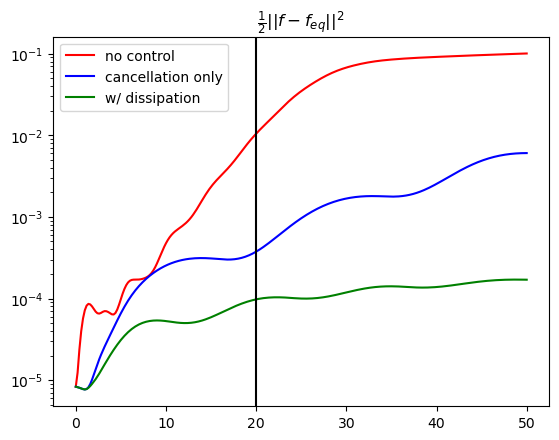}
    \subcaption{$\tau=1.0>\tau_{\mathrm{ref}}=0.6$}
    \end{subfigure}
    \quad
    \begin{subfigure}{0.3\textwidth}
    \includegraphics[width=\textwidth]{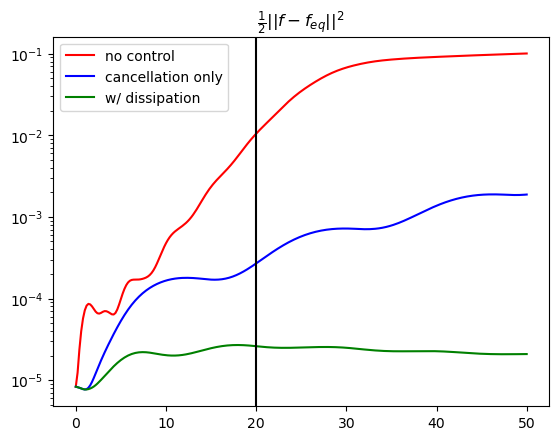}
    \subcaption{$\tau=0.8>\tau_{\mathrm{ref}}=0.6$}
    \end{subfigure}
    \quad
    \begin{subfigure}{0.3\textwidth}
    \includegraphics[width=\textwidth]{VM1D2V_ts_jdbde_perturb_tau0.6_lag0.6.png}
    \subcaption{$\tau=\tau_{\mathrm{ref}}=0.6$}
    \end{subfigure}
    \\
    \begin{subfigure}{0.3\textwidth}
    \includegraphics[width=\textwidth]{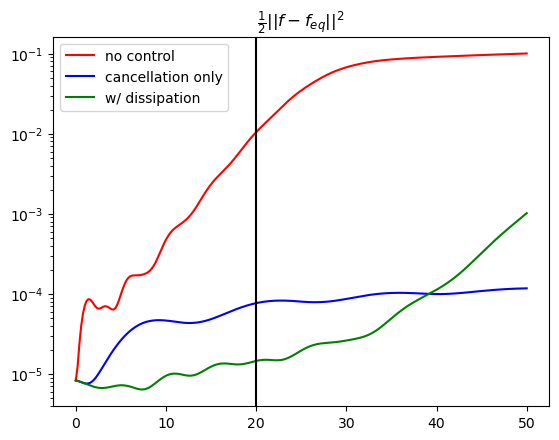}
    \subcaption{$\tau=0.4<\tau_{\mathrm{ref}}=0.6$}
    \end{subfigure}
    \quad
    \begin{subfigure}{0.3\textwidth}
    \includegraphics[width=\textwidth]{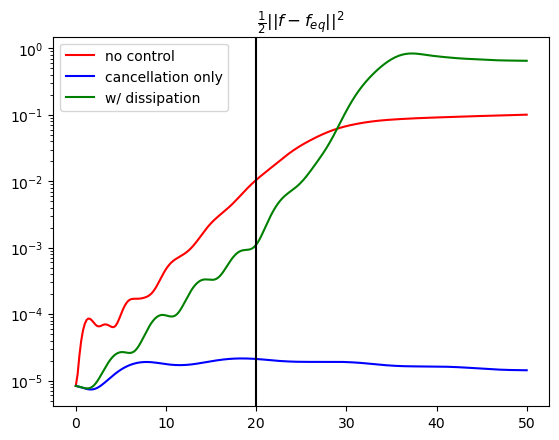}
    \subcaption{$\tau=0.2<\tau_{\mathrm{ref}}=0.6$}
    \end{subfigure}
    \caption{Two-stream instability. History of the $L^2$ perturbation under varying feedback delays $\tau$. No external field (red) vs.\ delayed cancellation control (blue) vs.\ cancellation--dissipation control (green). The dissipation networks \eqref{eq:dissipation} are trained at $\tau_{\mathrm{ref}}=0.6$ and then applied directly at different values of $\tau$ without correction. The vertical black line marks the terminal time of optimization.}
    \label{fig:two stream biased time lag}
\end{figure}

\begin{figure}[h!]
    \centering
    \begin{subfigure}{0.35\textwidth}
    \includegraphics[width=\textwidth]{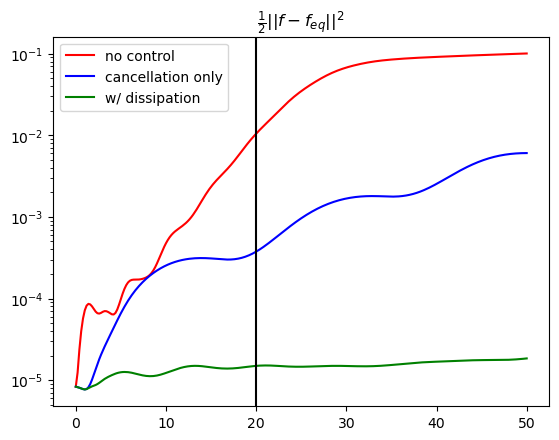}
    \subcaption{$\tau=1.0>\tau_{\mathrm{ref}}=0.6$}
    \end{subfigure}
    \quad
    \begin{subfigure}{0.35\textwidth}
    \includegraphics[width=\textwidth]{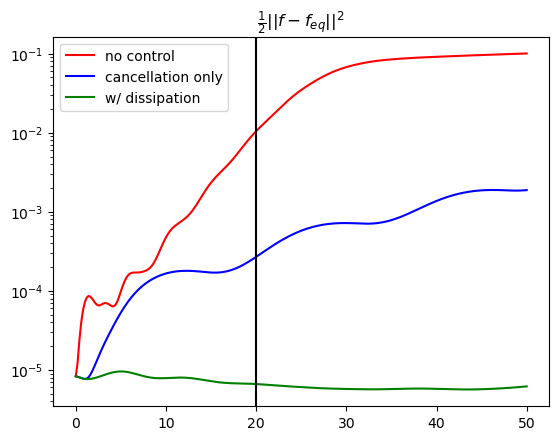}
    \subcaption{$\tau=0.8>\tau_{\mathrm{ref}}=0.6$}
    \end{subfigure}
    \\
    \begin{subfigure}{0.35\textwidth}
    \includegraphics[width=\textwidth]{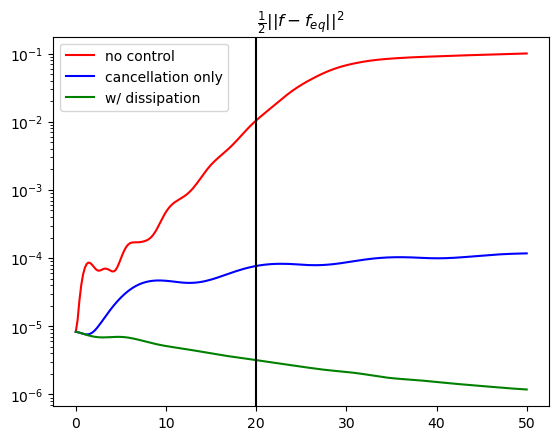}
    \subcaption{$\tau=0.4<\tau_{\mathrm{ref}}=0.6$}
    \end{subfigure}
    \quad
    \begin{subfigure}{0.35\textwidth}
    \includegraphics[width=\textwidth]{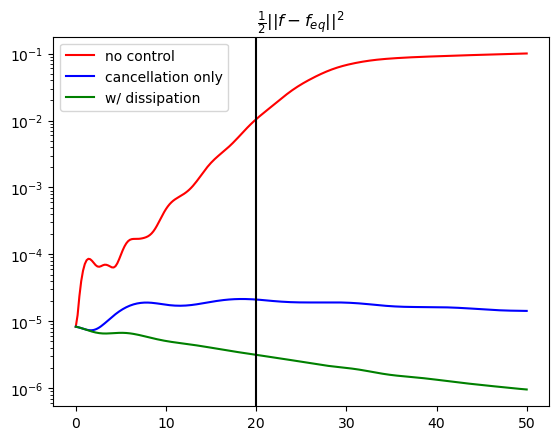}
    \subcaption{$\tau=0.2<\tau_{\mathrm{ref}}=0.6$}
    \end{subfigure}
    \caption{Two-stream instability. History of the $L^2$ perturbation under varying feedback delays $\tau$. No external field (red) vs.\ delayed cancellation control (blue) vs.\ cancellation--dissipation control (green). The dissipation networks \eqref{eq:dissipation} are trained at $\tau_{\mathrm{ref}}=0.6$, and their outputs are rescaled by $\tau/\tau_{\mathrm{ref}}$ in each test case.}
    \label{fig:two stream biased time lag corrected}
\end{figure}

In applications, the actual feedback delay $\tau$ may differ from the
reference value $\tau_{\mathrm{ref}}$ used during training. To assess the
robustness of the controller to such mismatch, we optimize
\eqref{eq:feedback operator} at $\tau_{\mathrm{ref}}=0.6$ and then apply
the learned controller at different values of $\tau$. As shown in
Figure~\ref{fig:two stream biased time lag}, when
$\tau>\tau_{\mathrm{ref}}$, the learned dissipation underestimates the
delay-induced error and therefore provides insufficient compensation.
The more delicate case occurs when $\tau<\tau_{\mathrm{ref}}$, for which
the learned dissipation can overcompensate and degrade stability.
A simple correction is to rescale the output of the dissipation networks
by $\tau/\tau_{\mathrm{ref}}$, without retraining. This scaling is
consistent with the first-order approximation of the time-integrated
cancellation error in \eqref{eq:cancellation err}, which is proportional
to the delay. Figure~\ref{fig:two stream biased time lag corrected}
shows that this correction is effective for both
$\tau<\tau_{\mathrm{ref}}$ and $\tau>\tau_{\mathrm{ref}}$.

\begin{figure}[h!]
    \centering
    \begin{subfigure}{0.3\textwidth}
    \includegraphics[width=\textwidth]{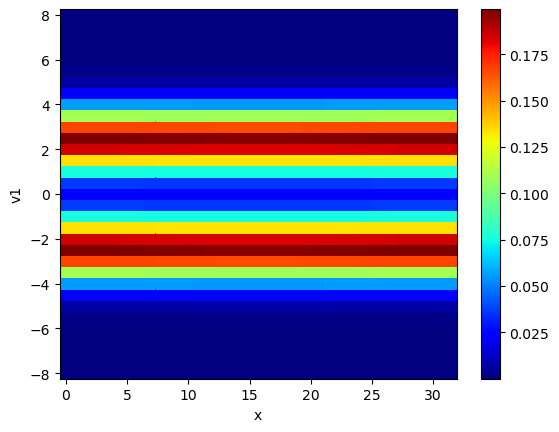}
    \subcaption{$t=0$, no control}
    \end{subfigure}
    \quad
    \begin{subfigure}{0.3\textwidth}
    \includegraphics[width=\textwidth]{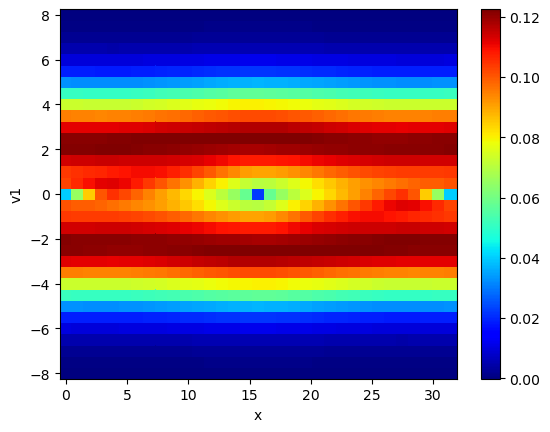}
    \subcaption{$t=50$, no control}
    \end{subfigure}
    \quad
    \begin{subfigure}{0.3\textwidth}
    \includegraphics[width=\textwidth]{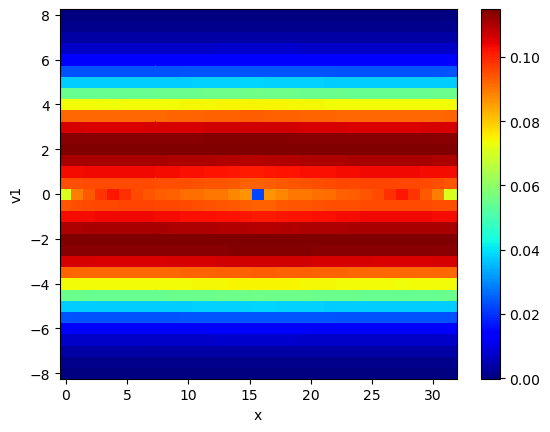}
    \subcaption{$t=100$, no control}
    \end{subfigure}
    \\
    \begin{subfigure}{0.3\textwidth}
    \includegraphics[width=\textwidth]{VM1D2V_ts_jdbde_f0.png}
    \subcaption{$t=0$, delayed cancellation}
    \end{subfigure}
    \quad
    \begin{subfigure}{0.3\textwidth}
    \includegraphics[width=\textwidth]{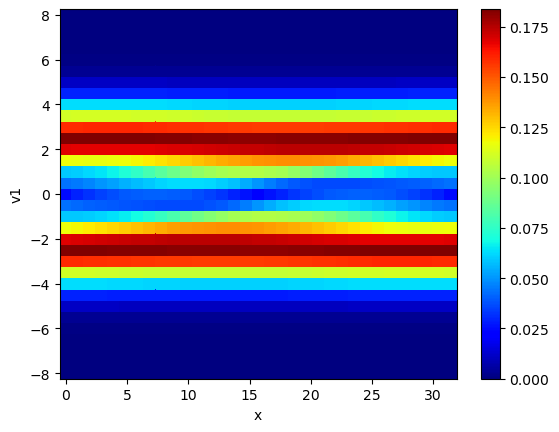}
    \subcaption{$t=50$, delayed cancellation}
    \end{subfigure}
    \quad
    \begin{subfigure}{0.3\textwidth}
    \includegraphics[width=\textwidth]{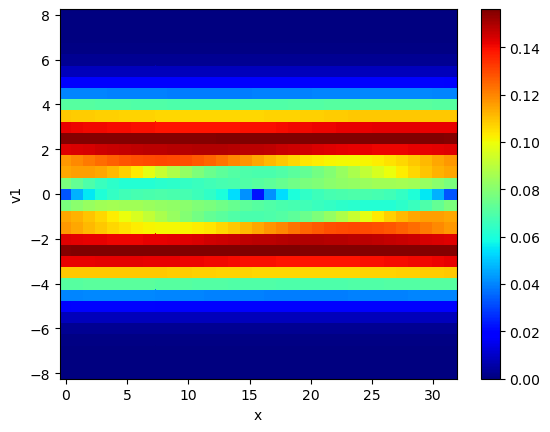}
    \subcaption{$t=100$, delayed cancellation}
    \end{subfigure}
    \\
    \begin{subfigure}{0.3\textwidth}
    \includegraphics[width=\textwidth]{VM1D2V_ts_jdbde_f0.png}
    \subcaption{$t=0$, operator control}
    \end{subfigure}
    \quad
    \begin{subfigure}{0.3\textwidth}
    \includegraphics[width=\textwidth]{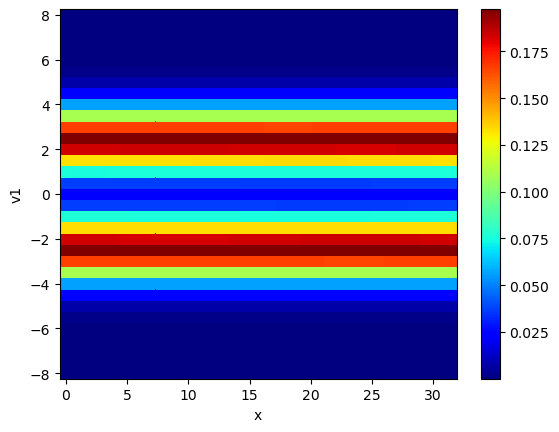}
    \subcaption{$t=50$, operator control}
    \end{subfigure}
    \quad
    \begin{subfigure}{0.3\textwidth}
    \includegraphics[width=\textwidth]{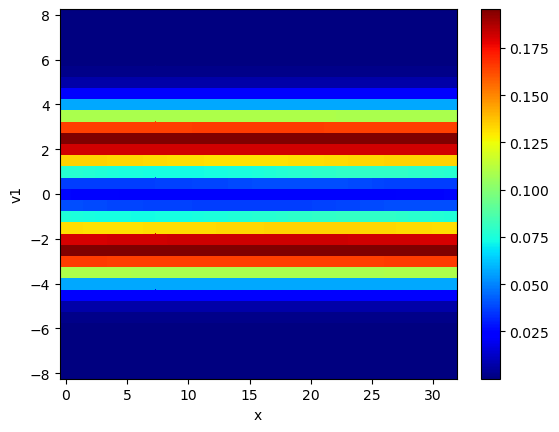}
    \subcaption{$t=100$, operator control}
    \end{subfigure}
    \caption{Two-stream instability. Evolution of the marginal distribution $\widetilde f(x,v_1,t)=\int f(x,v_1,v_2,t)\,dv_2$ for feedback delay $\tau=1.0$. No external field (first row) vs.\ delayed cancellation (second row) vs.\ cancellation--dissipation control (third row). The dissipation networks are trained at $\tau_{\mathrm{ref}}=0.6$ over $t\in[0,20]$, consistent with Figures~\ref{fig:two stream time lag perturb}--\ref{fig:two stream biased time lag corrected}, and their outputs are rescaled by $\tau/\tau_{\mathrm{ref}}=5/3$.}
    \label{fig:two stream time lag f}
\end{figure}

We further display the marginal distribution
\[
\widetilde f(x,v_1,t)
=
\int f(x,v_1,v_2,t)\,dv_2
\]
in Figure~\ref{fig:two stream time lag f} for the larger delay
$\tau=1.0$. Delayed cancellation alone provides noticeably weaker
suppression of the kinetic deformation, whereas the operator-based
controller better preserves the two-stream structure over long times.
Together, these results indicate that the learned correction can
compensate for a non-negligible feedback delay and that moderate delay
mismatch can be handled by a simple rescaling of the learned dissipation.

\subsection{Weibel-type current-filamentation instability}

We next consider the equilibrium
\[
\begin{split}
&f_{\mathrm{eq}}(v_1,v_2)
=
\frac{1}{2\pi\beta}
\exp\left(
-\frac{v_1^2}{\beta}
-\frac{(v_2-\bar v_2)^2}{\beta}
\right)
+
\frac{1}{2\pi\beta}
\exp\left(
-\frac{v_1^2}{\beta}
-\frac{(v_2+\bar v_2)^2}{\beta}
\right),\\
&E_{eq,1}=0,\qquad
E_{eq,2}=0,\qquad
B_{eq,3}=0,
\end{split}
\]
with $\beta=0.04$ and $\bar v_2=0.35$. This distribution consists of
two equal counterstreaming populations in the $v_2$ direction. Since the
spatial variation is along $x$, perpendicular to the streaming direction,
the resulting transverse electromagnetic instability is more precisely
described as a current-filamentation instability, closely related to the
Weibel instability. We therefore refer to this test as a
Weibel-type current-filamentation instability.

\begin{figure}[h!]
    \centering
    \begin{subfigure}{0.4\textwidth}
    \includegraphics[width=\textwidth]{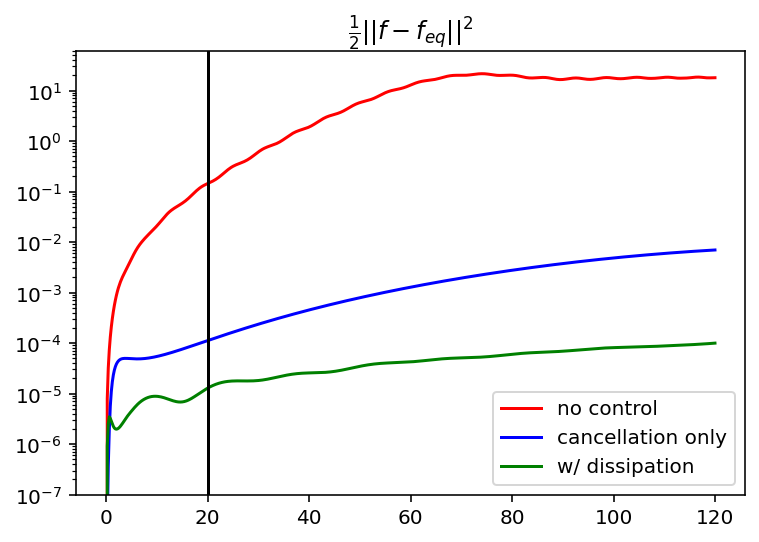}
    \subcaption{$L^2$ perturbation}
    \end{subfigure}
    \quad
    \begin{subfigure}{0.4\textwidth}
    \includegraphics[width=\textwidth]{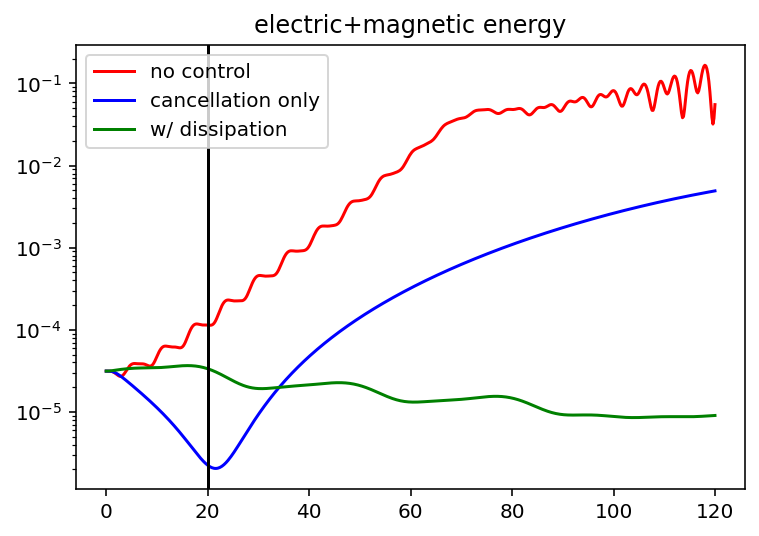}
    \subcaption{Electromagnetic energy}
    \end{subfigure}
    \caption{Weibel-type current-filamentation instability. Histories of the $L^2$ perturbation and electromagnetic energy over $t\in[0,120]$. No external field (red) vs.\ delayed cancellation control (blue) vs.\ cancellation--dissipation control (green). The feedback delay is $\tau=\tau_{\mathrm{ref}}=0.6$. The vertical black line marks the terminal time of optimization.}
    \label{fig:weibel time lag perturb}
\end{figure}

The computational domain is
\[
(x,v_1,v_2)\in[0,10\pi]\times[-1,1]\times[-1,1].
\]
The learning rate is set to $10^{-4}$ during the \texttt{Adagrad}
preprocessing stage and to $2\times10^{-5}$ during the subsequent
\texttt{Adam} iterations. The instability is seeded by the small initial
magnetic perturbation
\[
B_3(x,0)
=
0.002\cos\left(\frac{x}{5}\right).
\]
The anisotropic counterstreaming distribution responds to this magnetic
perturbation by developing transverse current structures, which reinforce
the magnetic perturbation and lead to rapid growth of the kinetic
deviation, as shown in Figure~\ref{fig:weibel time lag perturb}.

\begin{figure}[h!]
    \centering
    \begin{subfigure}{0.3\textwidth}
    \includegraphics[width=\textwidth]{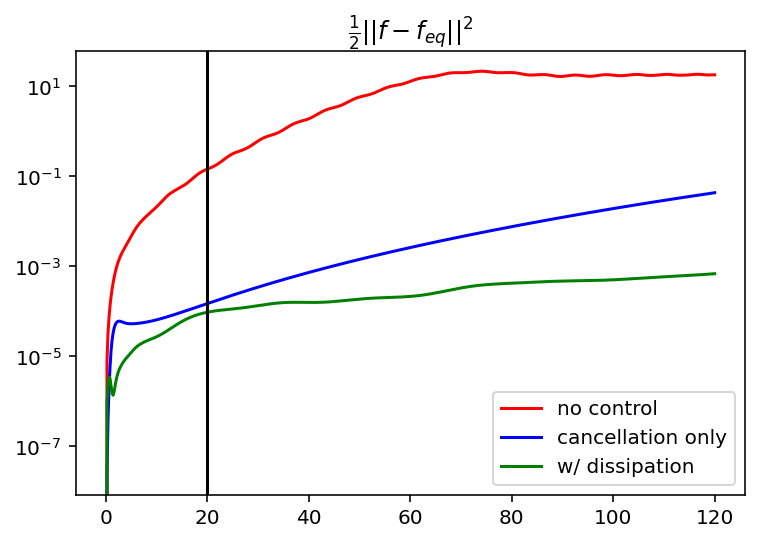}
    \subcaption{$\tau=0.8$, uncorrected dissip.}
    \end{subfigure}
    \quad
    \begin{subfigure}{0.3\textwidth}
    \includegraphics[width=\textwidth]{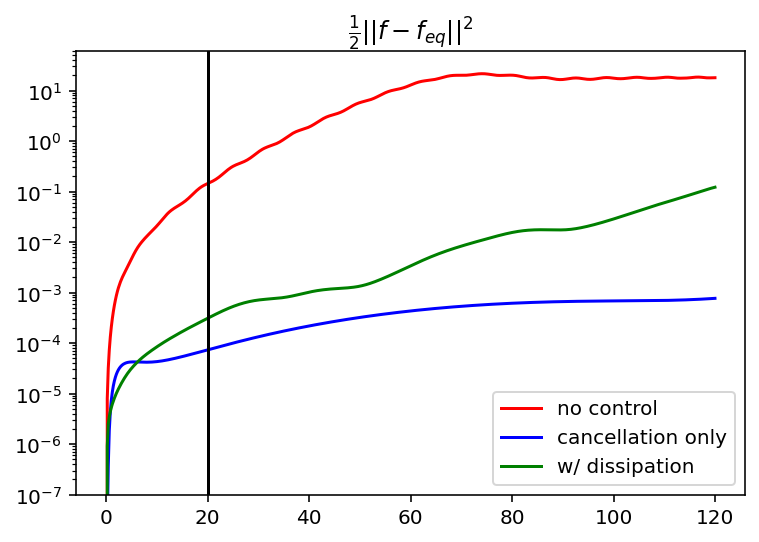}
    \subcaption{$\tau=0.4$, uncorrected dissip.}
    \end{subfigure}
    \quad
    \begin{subfigure}{0.3\textwidth}
    \includegraphics[width=\textwidth]{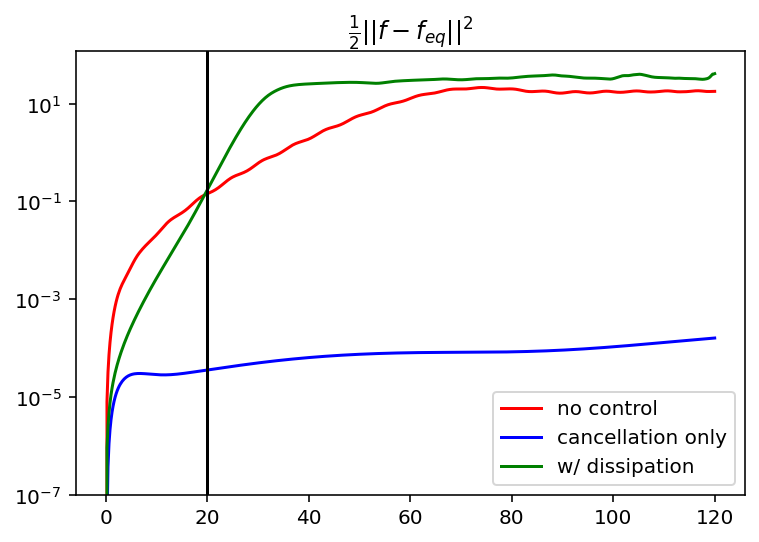}
    \subcaption{$\tau=0.2$, uncorrected dissip.}
    \end{subfigure}
    \\
    \begin{subfigure}{0.3\textwidth}
    \includegraphics[width=\textwidth]{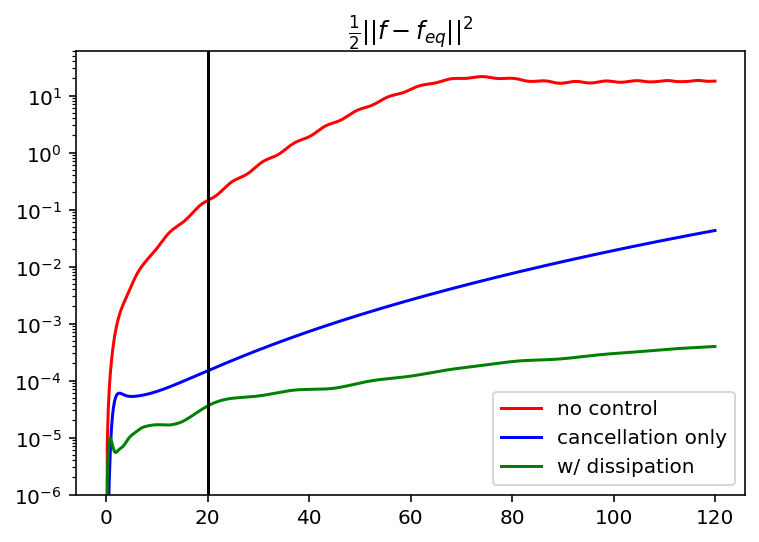}
    \subcaption{$\tau=0.8$, corrected dissip.}
    \end{subfigure}
    \quad
    \begin{subfigure}{0.3\textwidth}
    \includegraphics[width=\textwidth]{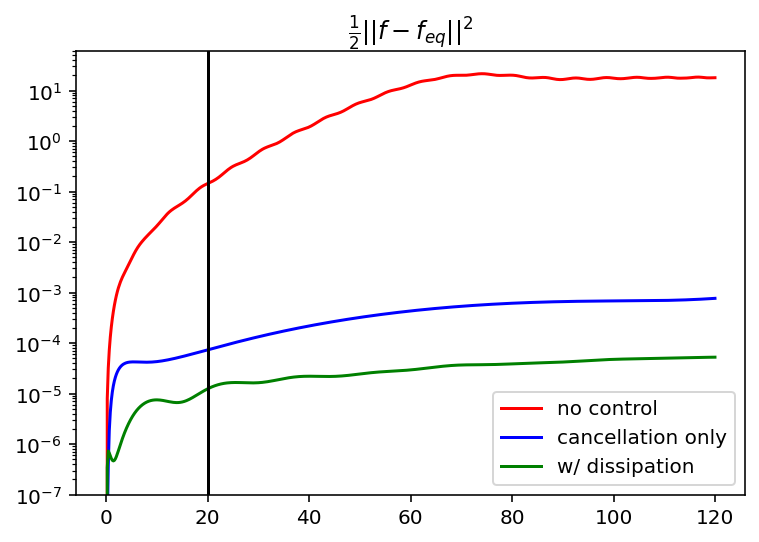}
    \subcaption{$\tau=0.4$, corrected dissip.}
    \end{subfigure}
    \quad
    \begin{subfigure}{0.3\textwidth}
    \includegraphics[width=\textwidth]{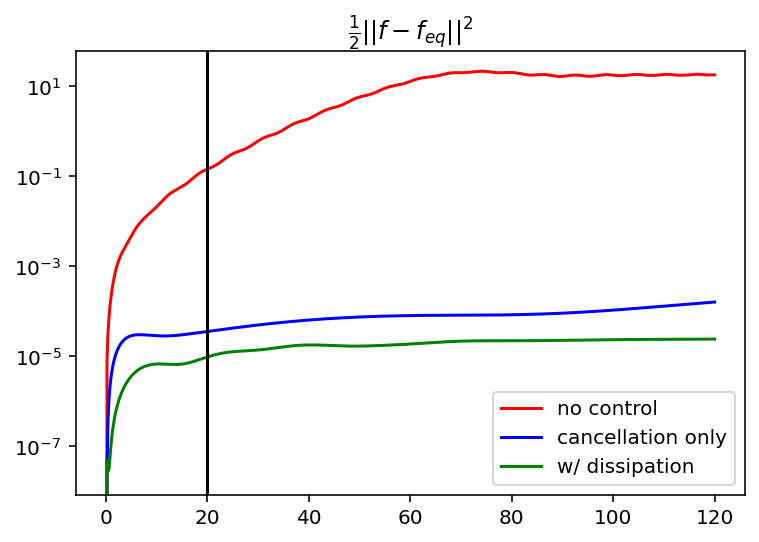}
    \subcaption{$\tau=0.2$, corrected dissip.}
    \end{subfigure}
    \caption{Weibel-type current-filamentation instability. History of the $L^2$ perturbation under varying feedback delays $\tau$. No external field (red) vs.\ delayed cancellation control (blue) vs.\ cancellation--dissipation control (green). The dissipation networks \eqref{eq:dissipation} are trained at $\tau_{\mathrm{ref}}=0.6$. In the first row, the learned dissipation is applied directly at $\tau\neq\tau_{\mathrm{ref}}$ without correction. In the second row, the network outputs are rescaled by $\tau/\tau_{\mathrm{ref}}$. The vertical black line marks the terminal time of optimization.}
    \label{fig:weibel biased time lag}
\end{figure}

The feedback controller \eqref{eq:feedback operator} is trained with
reference delay $\tau_{\mathrm{ref}}=0.6$ over $t\in[0,20]$. As in the
two-stream test, the learned controller substantially improves upon
delayed field cancellation. Figure~\ref{fig:weibel biased time lag}
shows that when $\tau>\tau_{\mathrm{ref}}$, the learned dissipation is
too weak to fully compensate for the larger delay, whereas
$\tau<\tau_{\mathrm{ref}}$ can lead to overcompensation. In both cases,
rescaling the learned dissipation by $\tau/\tau_{\mathrm{ref}}$
substantially reduces the effect of delay mismatch.

\begin{figure}[h!]
    \centering
    \begin{subfigure}{0.3\textwidth}
    \includegraphics[width=\textwidth]{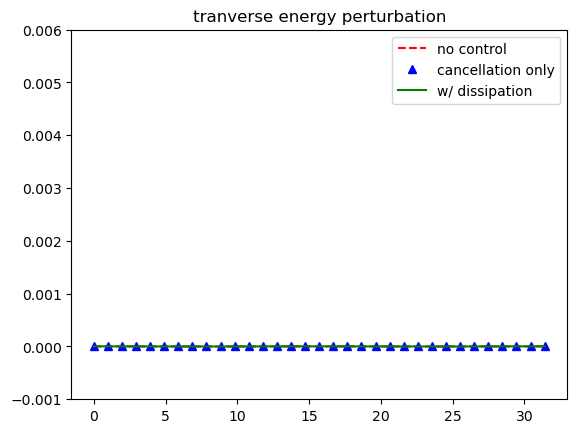}
    \subcaption{$t=0$}
    \end{subfigure}
    \quad
    \begin{subfigure}{0.3\textwidth}
    \includegraphics[width=\textwidth]{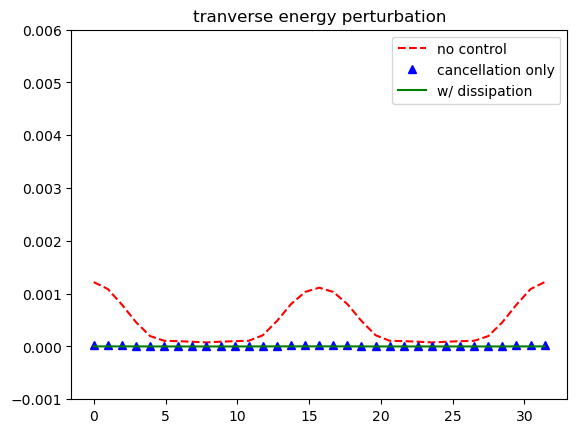}
    \subcaption{$t=60$}
    \end{subfigure}
    \quad
    \begin{subfigure}{0.3\textwidth}
    \includegraphics[width=\textwidth]{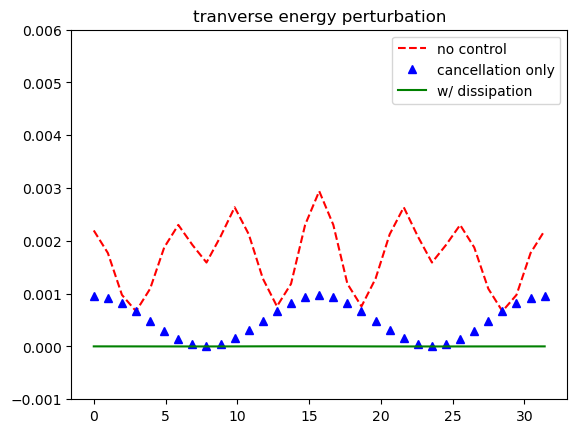}
    \subcaption{$t=120$}
    \end{subfigure}
    \caption{Weibel-type current-filamentation instability. Evolution of the transverse electromagnetic energy density $\frac{1}{2}(E_2(x,t)^2+B_3(x,t)^2)$ for feedback delay $\tau=0.8$. No external field (red) vs.\ delayed cancellation control (blue) vs.\ operator-based feedback control (green). The dissipation networks \eqref{eq:dissipation} are trained at $\tau_{\mathrm{ref}}=0.6$, with their outputs rescaled by $\tau/\tau_{\mathrm{ref}}=4/3$.}
    \label{fig:weibel transverse perturb}
\end{figure}

Figure~\ref{fig:weibel transverse perturb} shows the spatial profiles of
the local transverse electromagnetic energy density
\[
\frac{1}{2}\left(E_2(x,t)^2+B_3(x,t)^2\right)
\]
at selected times. Delayed field cancellation suppresses the transverse
field growth over a moderate time interval but gradually loses
effectiveness. In contrast, the cancellation--dissipation controller
continues to suppress the transverse perturbation well beyond the
optimization horizon.

\newpage

\vspace{0.5in}

\subsection{Effect of feedback-observation design}

To assess the role of observation design in long-time stabilization and
delay compensation, we compare the proposed gradient-based feedback
variables
\[
\left(
\delta J_1,\delta J_2,
\partial_x\delta E_2,\partial_x\delta B_3
\right)
\]
with the component-based alternative
\[
\left(
\delta J_1,\delta J_2,\delta E_2,\delta B_3
\right).
\]
The component-based choice provides a natural baseline because the
electromagnetic field components directly characterize the instantaneous field
state. However, the delay-compensation identities
\eqref{eq:cancellation err} show more specifically that the cancellation
errors are generated by $\delta J_1$, $\delta J_2$,
$\partial_x\delta E_2$, and $\partial_x\delta B_3$. These quantities
therefore encode the instantaneous field evolution responsible for the
delay-induced residual and are expected to provide more informative
feedback variables. In the following comparison, all network architectures
and training settings are kept unchanged; only the last two feedback
variables are replaced.

\begin{figure}[h!]
    \centering
    \begin{subfigure}{0.43\textwidth}
    \includegraphics[width=\textwidth]{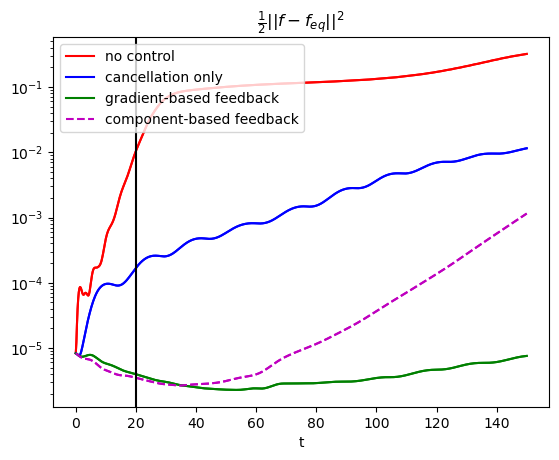}
    \subcaption{Two-stream: $L^2$ perturbation}
    \end{subfigure}
    \quad
    \begin{subfigure}{0.43\textwidth}
    \includegraphics[width=\textwidth]{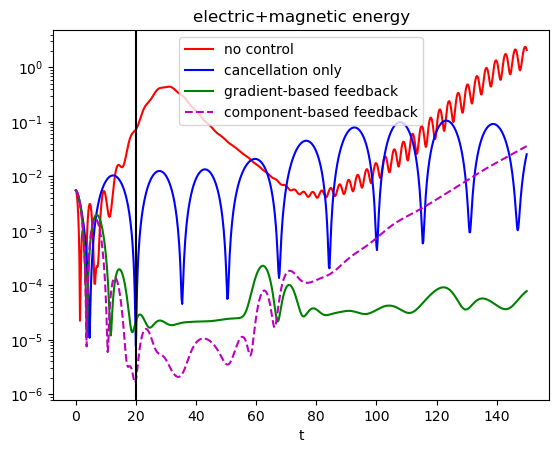}
    \subcaption{Two-stream: electromagnetic energy}
    \end{subfigure}
    \\
    \begin{subfigure}{0.43\textwidth}
    \includegraphics[width=\textwidth]{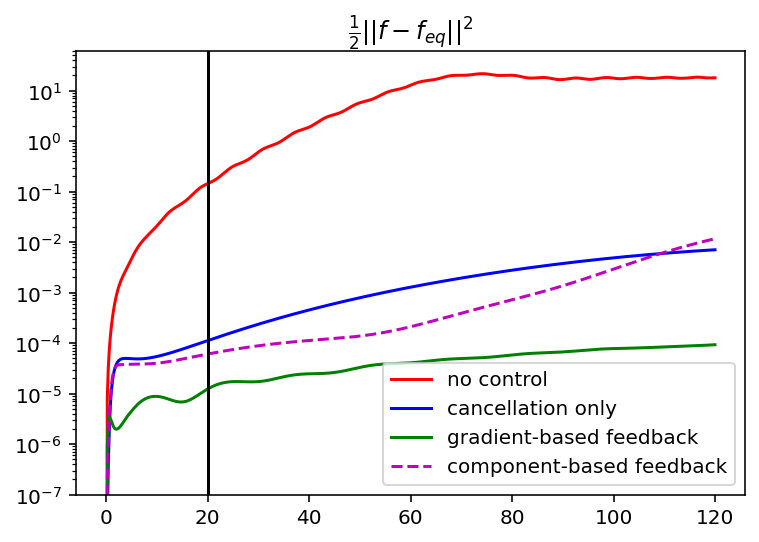}
    \subcaption{Weibel-type filamentation: $L^2$ perturbation}
    \end{subfigure}
    \quad
    \begin{subfigure}{0.43\textwidth}
    \includegraphics[width=\textwidth]{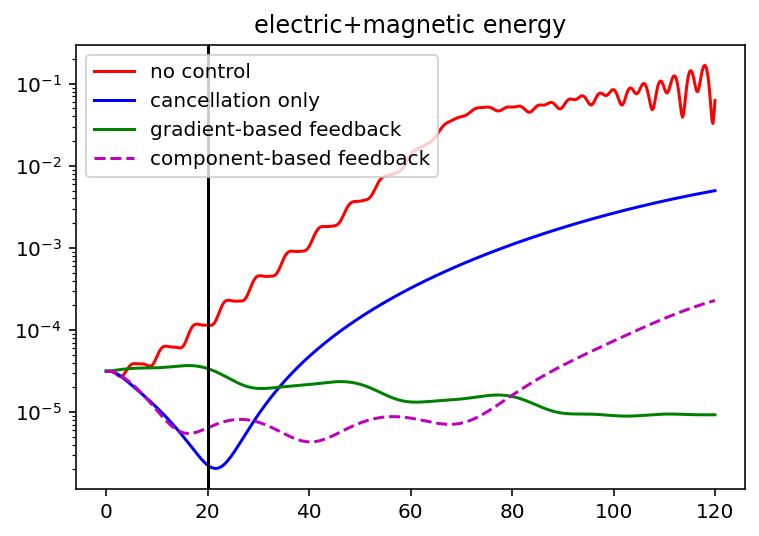}
    \subcaption{Weibel-type filamentation: electromagnetic energy}
    \end{subfigure}
    \caption{Comparison of gradient-based and component-based feedback for the two-stream instability (first row) and the Weibel-type current-filamentation instability (second row). Histories of the $L^2$ perturbation and electromagnetic energy over $t\in[0,150]$ for the two-stream instability and $t\in[0,120]$ for the Weibel-type instability. No external field (red) vs.\ delayed cancellation control (blue) vs.\ gradient-based delay compensation (green) vs.\ component-based delay compensation (pink). Both learned controllers use the operator structure \eqref{eq:feedback operator}--\eqref{eq:dissipation}; the component-based controller replaces $\partial_x\delta E_2$ and $\partial_x\delta B_3$ by $\delta E_2$ and $\delta B_3$, respectively. The feedback delay is $\tau=\tau_{\mathrm{ref}}=0.6$.}
    \label{fig:grad vs comp}
\end{figure}

Figure~\ref{fig:grad vs comp} presents the comparison for the two-stream
and Weibel-type current-filamentation instabilities. In both cases, the
component-based controller provides substantial short-time suppression but
degrades significantly beyond the optimization horizon. In contrast, the
gradient-based controller maintains substantially stronger long-time
stabilization of both the kinetic perturbation and the electromagnetic
energy. The consistency of this behavior across the two instability
mechanisms supports the observation design derived directly from the
Vlasov--Maxwell dynamics.

We note that spatial differentiation can amplify measurement noise, and
its practical implementation depends on sensor resolution and the
procedure used to reconstruct derivatives. Incorporating noisy
observations and regularized derivative estimation into the training
procedure is left for future work.

\subsection{Linear versus nonlinear operator architectures}

The main experiments above use dissipation networks with a linear operator
structure. We now examine whether introducing nonlinearity into the
operator architecture changes the control performance. As a representative
nonlinear architecture, we use a Fourier neural operator (FNO) with two
Fourier layers, 15 retained modes, and width 9 in each intermediate layer.
The number of retained modes is consistent with the spectral resolution
used in the preceding LNO experiments, and the biases are set to zero.
With this configuration, the FNO- and LNO-based controllers have comparable
numbers of trainable parameters. We replace the linear dissipation operator
in \eqref{eq:dissipation} by the FNO while keeping the remainder of the
training procedure unchanged. To account for the different optimization
behavior of the nonlinear architecture, the learning rates in the
\texttt{Adagrad} and \texttt{Adam} stages are set to
$5\times10^{-3}$ and $10^{-3}$, respectively.

\begin{figure}[h!]
    \centering
    \begin{subfigure}{0.43\textwidth}
    \includegraphics[width=\textwidth]{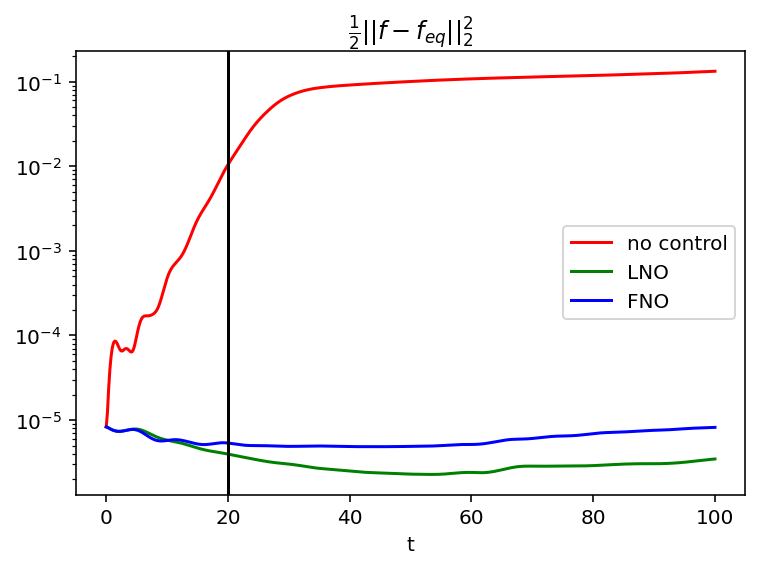}
    \subcaption{$L^2$ perturbation}
    \end{subfigure}
    \quad
    \begin{subfigure}{0.43\textwidth}
    \includegraphics[width=\textwidth]{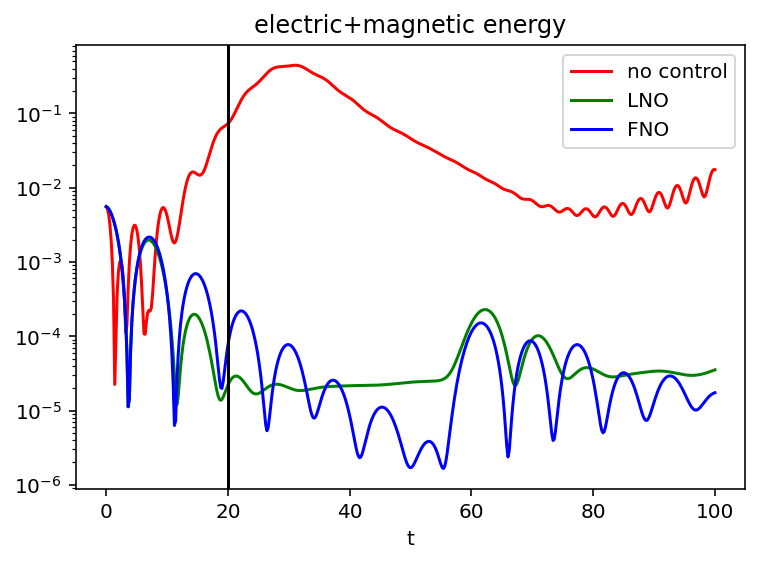}
    \subcaption{Electromagnetic energy}
    \end{subfigure}
    \caption{Comparison of LNO- and FNO-based controllers for the two-stream instability. Histories of the $L^2$ perturbation and electromagnetic energy over $t\in[0,100]$. No external field (red) vs.\ LNO-based control (green) vs.\ FNO-based control (blue). Both controllers are trained at $\tau_{\mathrm{ref}}=0.6$ over $t\in[0,20]$ and show comparable performance at the reference delay.}
    \label{fig:two stream lno vs fno}
\end{figure}

\begin{figure}[h!]
    \centering
    \begin{subfigure}{0.43\textwidth}
    \includegraphics[width=\textwidth]{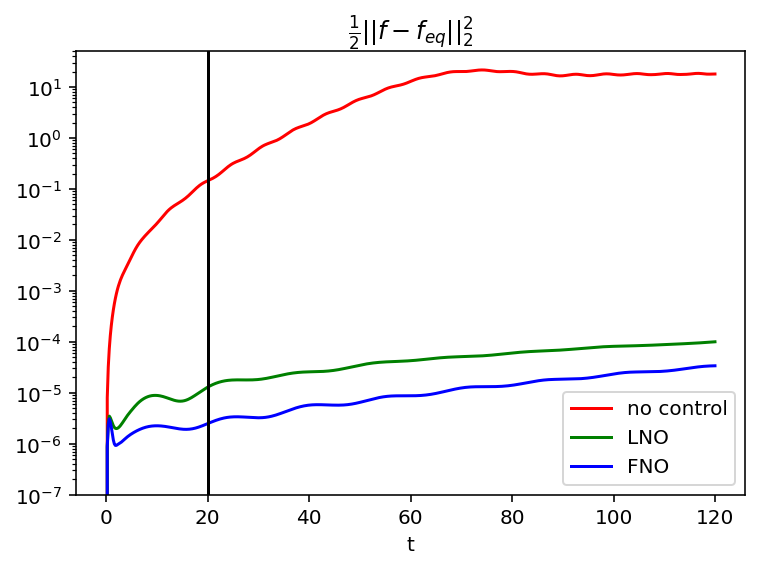}
    \subcaption{$L^2$ perturbation}
    \end{subfigure}
    \quad
    \begin{subfigure}{0.43\textwidth}
    \includegraphics[width=\textwidth]{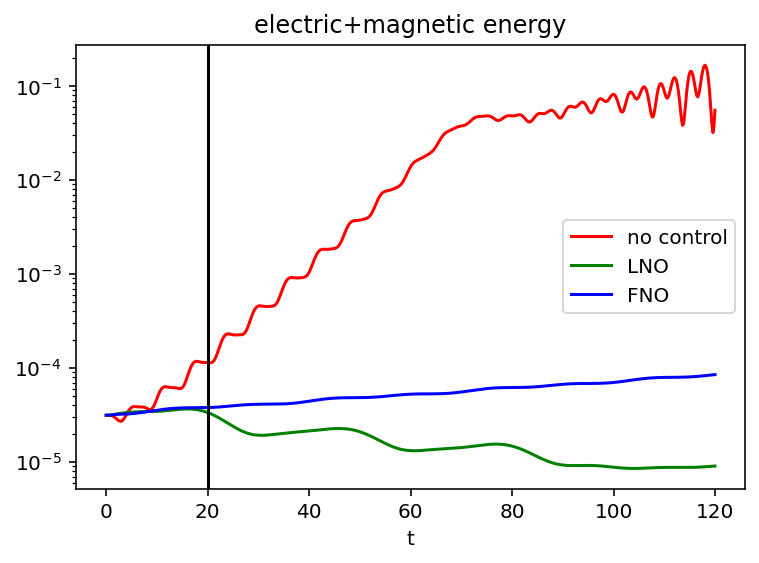}
    \subcaption{Electromagnetic energy}
    \end{subfigure}
    \caption{Comparison of LNO- and FNO-based controllers for the Weibel-type current-filamentation instability. Histories of the $L^2$ perturbation and electromagnetic energy over $t\in[0,120]$. No external field (red) vs.\ LNO-based control (green) vs.\ FNO-based control (blue). Both controllers are trained at $\tau_{\mathrm{ref}}=0.6$ over $t\in[0,20]$.}
    \label{fig:weibel lno vs fno}
\end{figure}

\begin{figure}[h!]
    \centering
    \begin{subfigure}{0.43\textwidth}
    \includegraphics[width=\textwidth]{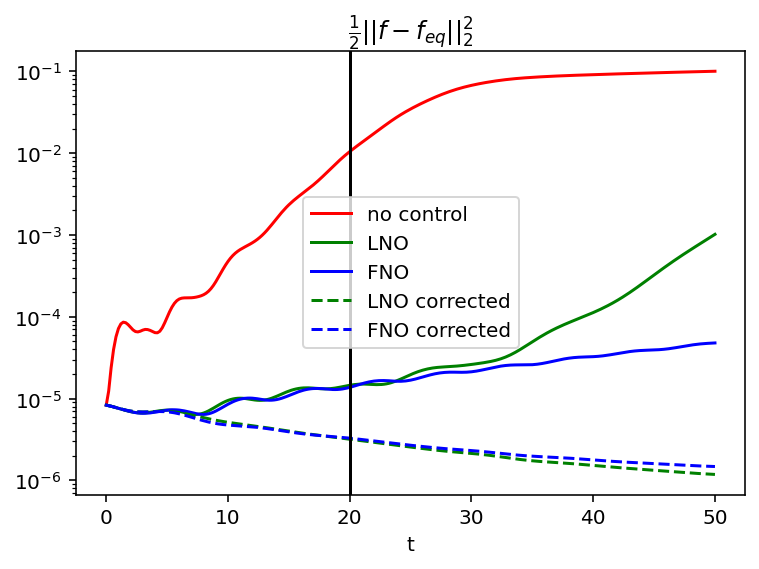}
    \subcaption{$\tau=0.4<\tau_{\mathrm{ref}}=0.6$}
    \end{subfigure}
    \quad
    \begin{subfigure}{0.43\textwidth}
    \includegraphics[width=\textwidth]{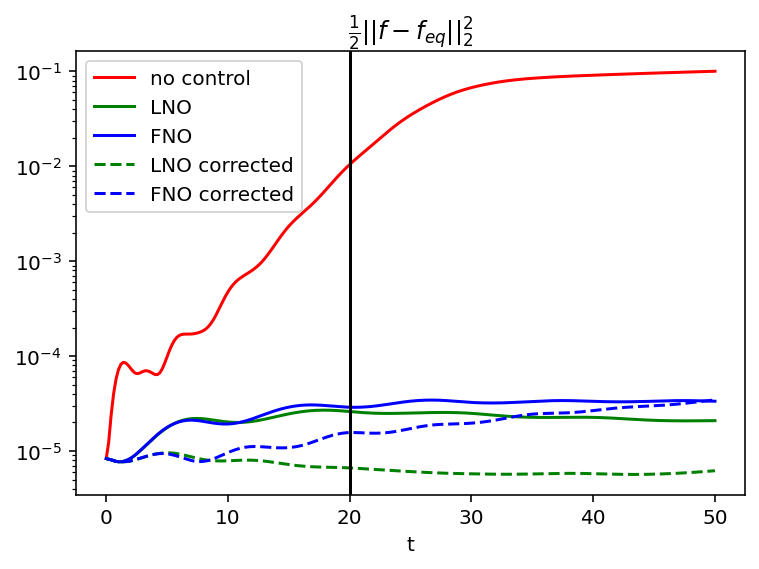}
    \subcaption{$\tau=0.8>\tau_{\mathrm{ref}}=0.6$}
    \end{subfigure}
    \caption{Comparison of LNO- and FNO-based controllers for the two-stream instability under varying feedback delays $\tau$. Histories of the $L^2$ perturbation for no external field (red), LNO-based control (green), FNO-based control (blue), corrected LNO-based control (green dashed), and corrected FNO-based control (blue dashed). The FNO is more robust for $\tau<\tau_{\mathrm{ref}}$, whereas the corrected LNO gives slightly stronger suppression in this test.}
    \label{fig:two stream lno vs fno biased time lag}
\end{figure}

\begin{figure}[h!]
    \centering
    \begin{subfigure}{0.43\textwidth}
    \includegraphics[width=\textwidth]{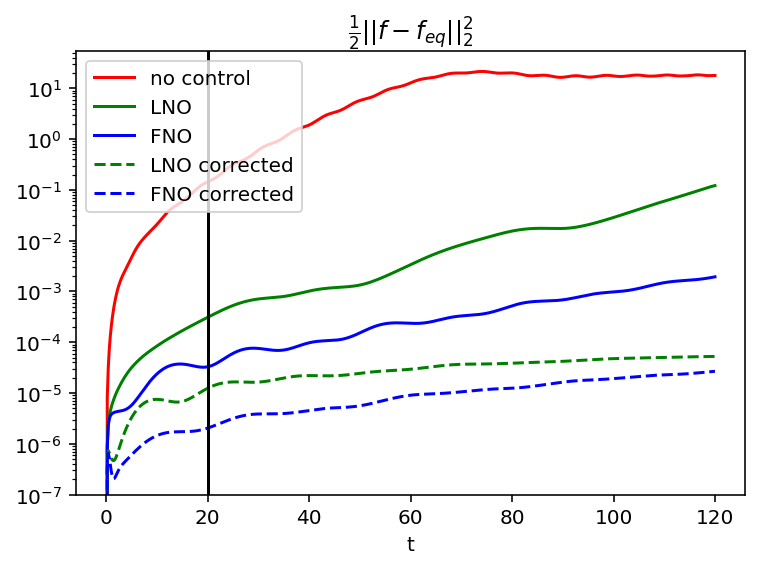}
    \subcaption{$\tau=0.4<\tau_{\mathrm{ref}}=0.6$}
    \end{subfigure}
    \quad
    \begin{subfigure}{0.43\textwidth}
    \includegraphics[width=\textwidth]{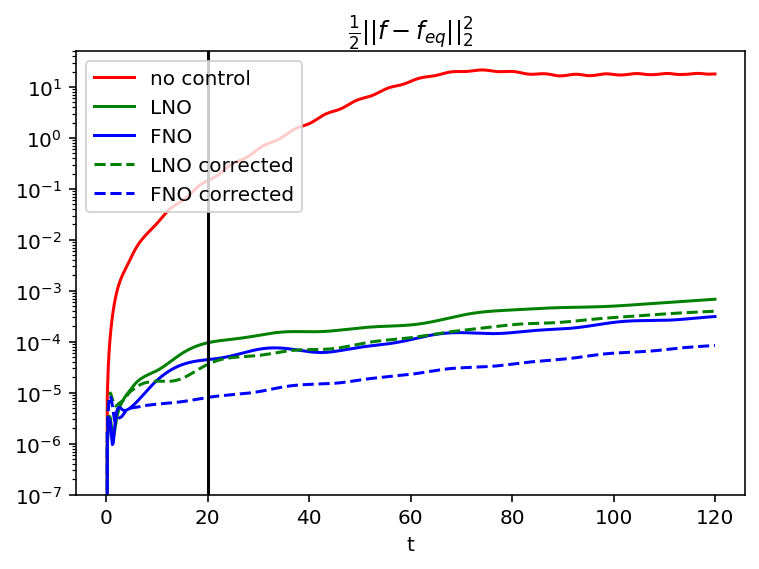}
    \subcaption{$\tau=0.8>\tau_{\mathrm{ref}}=0.6$}
    \end{subfigure}
    \caption{Comparison of LNO- and FNO-based controllers for the Weibel-type current-filamentation instability under varying feedback delays $\tau$. Histories of the $L^2$ perturbation for no external field (red), LNO-based control (green), FNO-based control (blue), corrected LNO-based control (green dashed), and corrected FNO-based control (blue dashed). In these tests, the FNO is less sensitive to delay mismatch.}
    \label{fig:weibel lno vs fno biased time lag}
\end{figure}

At the reference delay, Figures~\ref{fig:two stream lno vs fno} and
\ref{fig:weibel lno vs fno} show that the LNO and FNO provide comparable
stabilization. The difference becomes more pronounced when the actual
delay differs from the value used during training; see
Figures~\ref{fig:two stream lno vs fno biased time lag} and
\ref{fig:weibel lno vs fno biased time lag}. In particular, the FNO is
less sensitive in the tested cases when the actual delay is smaller than
the reference value, $\tau<\tau_{\mathrm{ref}}$, and avoids the pronounced
overcompensation observed for the LNO. The precise source of this
difference is not yet clear and may reflect both architectural and
optimization effects.

The delay-informed rescaling by $\tau/\tau_{\mathrm{ref}}$ improves both
architectures. After this correction, the relative performance is
problem-dependent: the LNO gives slightly stronger suppression in the
two-stream test, whereas the FNO remains moderately more effective in the
Weibel-type current-filamentation test. These results suggest a trade-off
between the structural simplicity of the LNO and the greater robustness
to delay mismatch observed for the FNO in the present experiments. The
LNO has a more interpretable structure and, through a flexible choice of
basis functions, can be adapted naturally to prescribed actuation modes
and nonuniform spatial discretizations. Without the delay-rescaling
correction, however, it is more sensitive to mismatch between the actual
delay and the value used during training. The standard FNO is less
sensitive to such mismatch in the present tests, although its Fourier
representation is most natural on regular grids.

\clearpage

\appendix

\section{Fourier and Laplace transform of linearized equation}\label{append:fourier laplace}

In this section, we derive the relation \eqref{eq:density perturb}. We denote
\[
\begin{split}
\widehat{\psi}(\xi,v_1,v_2,t) = \int\psi(x,v_1,v_2,t)e^{-i\xi x}dx, \qquad
\widehat{\phi}(\xi,t) = \int \phi(x,t) e^{-i\xi x}dx.
\end{split}
\]

Performing Fourier transform in $x$ in the linearized equation \eqref{eq: VM linearize} yields
\[
\begin{split}
\partial_t\widehat{\df}(\xi,v_1,v_2,t)+i\xi v_1\widehat{\df}(\xi,v_1,v_2,t) = &-\big(\hE_1(\xi,t)+\hH_1(\xi,t)\big) \partial_{v_1}\feq(v_1,v_2)\\
&-\big (\hB_3(\xi,t)+\hB_{ext}(\xi,t)\big)v_2 \partial_{v_1}\feq(v_1,v_2)\\
&-(\hE_2(\xi,t)+\hH_2(\xi,t))\partial_{v_2}\feq(v_1,v_2)+(\hB_3(\xi,t)+\hB_{ext}(\xi,t))v_1\partial_{v_2}\feq(v_1,v_2).
\end{split}
\]
Note that $\ds \int\partial_{v_2}\feq(\xi,v_2)dv_2 = 0$, integrating the above equation with respect to $v_2$ eliminates the last two terms on the right,
\[
\begin{split}
\partial_t\int \widehat{\df}(\xi,v_1,v_2,t)dv_2
+i\xi v_1\int \widehat{\df}(\xi,v_1,v_2,t)dv_2 &= 
-\big(\hE_1(\xi,t)+\hH_1(\xi,t)\big)\int \partial_{v_1}\feq (v_1,v_2)dv_2\\
&-\big(\hB_3(\xi,t)+\hB_{ext}(\xi,t)\big)\int v_2\partial_{v_1}\feq(v_1,v_2)dv_2.
\end{split}
\]
Integration in time leads to
\[
\begin{split}
\int\widehat{\df}(\xi,v_1,v_2,t)dv_2 &= e^{-i\xi v_1 t}\int\widehat{\df}(\xi,v_1,v_2,0)dv_2\\
&-\int^t_0 e^{-i\xi v_1(t-\tau)}\big(\hE_1(\xi,\tau)+\hH_1(\xi,\tau)\big)[\int \partial_{v_1} \feq(v_1,v_2)dv_2]d\tau\\
&-\int^t_0 e^{-i\xi v_1(t-\tau)}\big(\hB_3(\xi,\tau)+\hB_{ext}(\xi,\tau)\big)[\int v_2\partial_{v_1}\feq(v_1,v_2)dv_2]d\tau
\end{split}
\]
Further integrating the above equation with respect to $v_1$, we notice that
\[
\begin{split}
\int\widehat{\df}(\xi,v_1,v_2,t)dv_1dv_2 = \widehat{\drho}(\xi,t),
\end{split}
\]
and 
\[
\begin{split}
\iint e^{-i\xi v_1 t}\widehat{\df}(\xi,v_1,v_2,0)dv_2dv_1 & = \iiint \df_0(x,v_1,v_2)e^{-i\xi(x+v_1 t)}dv_1dv_2dx\\
& = \iiint \df_0(x-v_1 t,v_1,v_2)e^{-i\xi x}dv_1 dv_2 dx\\
& = \widehat{\drho}_{F}(\xi,t),
\end{split}
\]
with $\ds \drho_{F}(x,t) = \int\df_0(x-v_1 t,v_1,v_2)dv_1dv_2$ being the free-stream density perturbation. Thus,
\begin{equation}\label{eq:hrho}
\widehat{\drho}(\xi,t) = \widehat{\drho}_{F}(\xi,t)-i\xi\int^t_0\big(\hE_1(\xi,\tau)+\hH_1(\xi,\tau)\big)\widehat{U}(\xi,t-\tau)d\tau-i\xi\int^t_0\big(\hB_3(\xi,\tau)+\hB_{ext}(\xi,\tau)\big)\widehat{U_v}(\xi,t-\tau)d\tau,
\end{equation}
where $U$ and $U_v$ are defined through Fourier transforms,
\[
\begin{split}
\widehat{U}(\xi,t) &= t\int \feq(v_1,v_2)e^{-i\xi t v_1}dv_1 dv_2= t\int\widehat{\feq}(\xi t,v_2)dv_2,\\ 
\widehat{U_v}(\xi,t) &= t\int v_2\feq(v_1,v_2)e^{-i\xi t v_1}dv_1 dv_2= t\int v_2\widehat{\feq}(\xi t,v_2)dv_2.
\end{split}
\]

Finally, performing Laplace transform in time in \eqref{eq:hrho} yields
\begin{equation}\label{eq:L-F rho}
\begin{split}
\laplace[\hrho(\xi,\cdot)](s) &= \laplace[\widehat{\drho}_{F}(\xi,\cdot)]-i\xi\big(\laplace[\hE_1(\xi,\cdot)](s)+\laplace[\hH_1(\xi,\cdot)](s)\big)\laplace[\widehat{U}(\xi,\cdot)](s)\\
&\hspace{5em}-i\xi\big(\laplace[\hB_3(\xi,\cdot)](s)+\laplace[\hB_{ext}(\xi,\cdot)](s)\big)\laplace[\widehat{U_v}(\xi,\cdot)](s)\\
& = \laplace[\widehat{\drho}_{F}(\xi,\cdot)]-\laplace[\widehat{\drho}(\xi,\cdot)](s)\laplace[\widehat{U}(\xi,\cdot)](s)-i\xi\laplace[\hH_1(\xi,\cdot)](s)\laplace[\widehat{U}(\xi,\cdot)](s)\\
&\hspace{5em}-i\xi\big(\laplace[\hB_3(\xi,\cdot)](s)+\laplace[\hB_{ext}(\xi,\cdot)](s)\big)\laplace[\widehat{U_v}(\xi,\cdot)](s).
\end{split}
\end{equation}
The second line is obtained by using 
\[
\partial_x E_1 = \drho\leadsto i\xi\hE_1 = \widehat{\drho}.
\]
Then the desired relation \eqref{eq:density perturb} is obtained by rearranging \eqref{eq:L-F rho}.

{\bf Acknowledgment:} We would like to thank Prof. Lukas Einkemmer for insightful discussions on the 1D2V model and linear stability analysis. 


\bibliographystyle{plain}
\bibliography{Ref}
\end{document}